\documentclass[final]{siamltex}

\usepackage{amsmath,amssymb}
\usepackage{latexsym}
\usepackage{subfigure}
\usepackage{booktabs}
\usepackage{framed}
\usepackage{algorithm}
\usepackage{algorithmic}
\usepackage{graphicx}
\usepackage{multirow}

\newcommand{\rf}[1]{\kern 0cm{(\ref{#1}})}

\newcommand{\ct}[1]{\kern 0cm{\cite{#1}}}

\title{Reciprocal and Positive Real Balanced Truncations
	for Model Order Reduction of Descriptor Systems\thanks{This work
	was supported in part by JSPS KAKENHI Grant Number 16K00073.}}

\author{Yuichi Tanji\thanks{Faculty of Engineering and Design,
		Kagawa University, Takamatsu, Kagawa 761-0396,
		Japan ({\tt tanji@eng.kagawa-u.ac.jp}).}}

\begin{document}

\maketitle

\begin{abstract}
Model order reduction algorithms for large-scale descriptor systems
are proposed using balanced truncation, in which
symmetry or block skew symmetry (reciprocity) and the positive realness of
the original transfer matrix are preserved.
Two approaches based on standard and
generalized algebraic Riccati equations are proposed.
To accelerate the algorithms, a fast Riccati solver,
RADI (alternating directions implicit [ADI]-type
iteration for Riccati equations), is also introduced.
As a result, the proposed methods are general and efficient as
a model order reduction algorithm for descriptor systems
associated with electrical circuit networks.
\end{abstract}

\begin{keywords}
model order reduction, balanced truncation, Riccati equation, ADI
\end{keywords}

\begin{AMS}
15A21, 15A22, 15A23, 15A24
\end{AMS}

\pagestyle{myheadings}
\thispagestyle{plain}
\markboth{YUICHI TANJI}{Reciprocal and Positive Real Balanced Truncations}

%
%

\section{Introduction}

Consider the following linear time invariant system:
\begin{eqnarray}
E_0 \frac{dx(t)}{dt} &=& A_0 x(t) + B_0 u(t), \quad
z(t) = C_0 x(t),
\label{eqn:des1}
\end{eqnarray}
where
$E_0 \in \mathbb{R}^{n,n}, A_0 \in \mathbb{R}^{n,n}, B_0 \in \mathbb{R}^{n,m},
C_0 \in \mathbb{R}^{m,n}$,
$x(t) \in \mathbb{R}^{n}$ is the state, $u(t) \in \mathbb{R}^{m}$ is the input,
and $z(t) \in \mathbb{R}^{m}$ is the output.
This is referred to as a descriptor system in the control community.
The form of \rf{eqn:des1} is typical for linear passive
electrical circuits in which coefficient matrices $E_0$ and $A_0$ are symmetric
and transfer matrix $H(s) = C_0 \left( sE_0 - A_0 \right)^{-1} B_0$
is symmetric or block skew symmetric. The symmetric or
block skew symmetric property of the transfer matrix is referred to as
reciprocity in circuit theory.
The transfer matrix $H(s)$ is assumed to be positive real,
whose property is referred to as passivity.
Furthermore, matrix $E_0$ is assumed to be singular, which is also
typical for linear electrical networks.

Interconnect networks in integrated circuits, packages, and printed circuit boards
are mathematically described by a descriptor system. Then,
the system becomes extremely large scale, which effectively
prohibits simulation of the
interconnect networks combined with other circuit blocks, which is necessary for
design of the electronics system including the integrated circuits, packages,
and printed circuit boards.
Therefore, the descriptor system should be of a small size
without losing behavior from zero to a specified frequency.
Model order reduction (MOR) methods fulfill these requirements, and
Krylov subspace methods \ct{kry6}, \ct{kry4}, \ct{kry5},
\ct{pvl} are powerful MOR methods that can provide an
accurate reduced-order model for a problem; however, these methods
cannot guarantee stability when connected to other linear networks.
Several methods \ct{index}, \ct{prima} are also based on the Krylov subspace,
in which coefficient matrices $E_0$ and $A_0$ of \rf{eqn:des1}
are written in block skew symmetric form to guarantee
passivity; however, reciprocity is not guaranteed.
If the coefficient matrices are written in symmetric form,
reciprocity is guaranteed but passivity is not.
Therefore, passivity and reciprocity are not simultaneously satisfied
with Krylov subspace methods.

Positive real balanced truncation (PRBT) \ct{phil}, \ct{wong2}
is more accurate than Krylov subspace methods, especially at high frequencies.
PRBT is an extension of balanced truncation \ct{moore}
in which algebraic Riccati equations (AREs) are solved.
Passivity-preserving balanced truncation for electrical
circuits (PABTEC) \ct{stykel} provides a passive and reciprocal reduced-order model
in which generalized AREs (GAREs)
are solved by the Newton method with Cholesky factorized alternating direction
implicit iteration to solve Lyapunov equations \ct{benner3}.
Note that PABTEC is applicable
to index-1 and index-2 descriptor systems.
Although reciprocity is guaranteed by preserving the block skew
symmetry of the transfer matrix, symmetric cases are not considered in PABTEC.

Thus, in this paper, we present reciprocal and positive real balanced
truncations (RPRBTs), in which reciprocity is guaranteed
for symmetric and block skew symmetric cases for
index-1 and index-2 systems.
Moreover, a fast Riccati solver, i.e., RADI \ct{benner2}, is introduced.
RADI solves the following ARE:
\begin{eqnarray}
A^T X + X A + X B B^T X + C^T C = 0,
\label{eqn:are}
\end{eqnarray}
where $A \in \mathbb{R}^{p,p}, B \in \mathbb{R}^{p,m}$, and $C \in \mathbb{R}^{m,p}$.
RADI provides an ARE solution that is equivalent to solutions obtained by
three seemingly different methods \ct{subric}, \ct{simon}, \ct{wong2}.
In addition, RADI
is the fastest among these methods.
Furthermore, RADI can be applied to solving the following GARE:
\begin{eqnarray}
A^T X E + E^T X A + E^T X B B^T X E + C^T C = 0,
\label{eqn:gare}
\end{eqnarray}
where
$E \in \mathbb{R}^{n,n}, A \in \mathbb{R}^{n,n}, B \in \mathbb{R}^{n,m}$,
and $C \in \mathbb{R}^{m,n}$.
Thus, we present two methods based on an ARE or GARE that are applicable to
both index-1 and index-2 systems.

RADI provides an accurate ARE solution in which
Ritz values are used as the shift parameters \ct{benner}.
However, it has been suggested \ct{tanji} that Ritz values
are insufficient for low-frequency accurate PRBT with quadratic
ADI (QADI) \ct{wong2},
which is one of the three methods equivalent to RADI.
To improve accuracy at low frequencies, which is important for
electrical circuit simulations, a small negative
constant value is used together with complex Ritz values as a shift
in QADI.
Furthermore, we extend the discussion to provide useful shift selections
to generate low-frequency accurate reduced-order models using RADI.

The remainder of this paper is organized as follows.
In Section 2, we provide definitions of reciprocity
for index-1 and index-2 descriptor systems.
In Section 3, RPRBT for an index-1 system is provided based on the ARE.
In Section 4, RPRBT for index-1 and index-2 systems is provided based on the GARE,
and RPRBT for the index-2 system is also provided based on the ARE.
In Section 5, the shift selections of RADI are obtained to solve the ARE and GARE.
We give numerical examples in Section 6 to demonstrate
the effectiveness of the proposed method.
Finally, conclusions are presented in Section 7.

The following notations are used in this paper.
$A^{-1}$, $A^{T}$, and $A^{*}$ represent the inverse, transpose,
and conjugate transpose
of matrix $A$, respectively.

%
%

\section{Reciprocity}

Reciprocity is a fundamental principle of linear passive networks.
In passive reduced-order interconnect macromodeling algorithm (PRIMA) \ct{prima},
the coefficient matrices $E_0$ and $A_0$ are
written in block skew symmetric form to
guarantee the positive realness of the transfer function.
As positive realness is guaranteed by solving
the Lur'e equation or an ARE in PRBT,
we do not need to write the descriptor system in block skew
symmetric form.
Therefore, matrices $E_0$ and $A_0$ are written in symmetric form.
The coefficient matrices of \rf{eqn:des1} are expressed as follows
for the impedance matrix:
\begin{eqnarray}
\begin{aligned}
E_0 &=
\left[
\begin{array}{cc}
A_{\mathcal C} {\mathcal C} A_{\mathcal C}^T & 0 \\
0 & -{\mathcal L} \\
\end{array}
\right], \quad
A_0 =
\left[
\begin{array}{cc}
-A_{\mathcal G} {\mathcal G} A_{\mathcal G}^T & -A_{\mathcal L} \\
-A_{\mathcal L}^T & 0
\end{array}
\right], \\
B_0 &=
\left[
\begin{array}{c}
-A_{\mathcal I} \\
0
\end{array}
\right], \quad
C_0 = B_0^T,
\end{aligned}
\label{eqn:z}
\end{eqnarray}
as follows for the admittance matrix:
\begin{eqnarray}
\begin{aligned}
E_0 &=
\left[
\begin{array}{ccc}
A_{\mathcal C} {\mathcal C} A_{\mathcal C}^T & 0 & 0 \\
0 & -{\mathcal L} & 0\\
0 & 0 & 0\\
\end{array}
\right], \quad
A_0 =
\left[
\begin{array}{ccc}
-A_{\mathcal G} {\mathcal G} A_{\mathcal G}^T & -A_{\mathcal L} & -A_{\mathcal V} \\
-A_{\mathcal L}^T & 0 & 0 \\
-A_{\mathcal V}^T & 0 & 0 \\
\end{array}
\right], \\
B_0 &=
\left[
\begin{array}{c}
0 \\
0 \\
I_m \\
\end{array}
\right], \quad
C_0 =
\left[
\begin{array}{ccc}
0 & 0 & -I_m \\
\end{array}
\right],
\end{aligned}
\label{eqn:y}
\end{eqnarray}
and as follows for the hybrid matrix:
\begin{eqnarray}
\begin{aligned}
E_0 &=
\left[
\begin{array}{ccc}
A_{\mathcal C} {\mathcal C} A_{\mathcal C}^T & 0 & 0 \\
0 & -{\mathcal L} & 0\\
0 & 0 & 0\\
\end{array}
\right], \quad
A_0 =
\left[
\begin{array}{ccc}
-A_{\mathcal G} {\mathcal G} A_{\mathcal G}^T & -A_{\mathcal L} & -A_{\mathcal V} \\
-A_{\mathcal L}^T & 0 & 0 \\
-A_{\mathcal V}^T & 0 & 0 \\
\end{array}
\right], \\
B_0 &=
\left[
\begin{array}{cc}
-A_{\mathcal I} & 0 \\
0 & 0 \\
0 & I_{m/2} \\
\end{array}
\right], \quad
C_0 =
\left[
\begin{array}{ccc}
-A_{\mathcal I} & 0 & 0 \\
0 & 0 & -I_{m/2} \\
\end{array}
\right].
\end{aligned}
\label{eqn:h}
\end{eqnarray}
In \rf{eqn:z}-\rf{eqn:h}, ${\mathcal G}$, ${\mathcal L}$, and
${\mathcal C}$ are the conductance, inductance, and capacitance matrices,
respectively. In addition, $A_{\mathcal G}$, $A_{\mathcal C}$, $A_{\mathcal L}$,
$A_{\mathcal I}$, and $A_{\mathcal V}$ are
the incidence matrices to conductors, capacitors,
inductors, and independent current and voltage sources, respectively.
Here, $I_m$ is an $m \times m$ identity matrix.
Then, linear passive RLC networks satisfy the following theorem.
\begin{theorem}
\label{th:recipro}
The linear passive networks expressed by \rf{eqn:des1} are reciprocal.
\end{theorem} \\
{\it Proof.}~It is trivial that
the transfer matrix associated with impedance and admittance
are symmetric.
Thus, the hybrid matrix is proved to be block
skew symmetric.
The hybrid matrix $H(s)$ is given as follows:
\begin{eqnarray}
\begin{aligned}
H(s) &=
I^\circ_m
B_0^T
\left(
s E_0 - A_0 \right)^{-1}
B_0
=
\left[
\begin{array}{cc}
H_{11} & H_{12} \\
H_{21} & H_{22}
\end{array}
\right],
\end{aligned}
\label{eqn:gh}
\end{eqnarray}
where
\begin{eqnarray*}
I^\circ_m &=&
\left[
\begin{array}{cc}
I_{m/2} & 0 \\
0 & -I_{m/2}
\end{array}
\right].
\end{eqnarray*}
As $B_0^T \left( s E_0 - A_0 \right)^{-1} B_0$ is symmetric,
$H_{11} = H_{11}^T$, $H_{22} = H_{22}^T$, and $H_{12} = - H_{21}^T$. \\
\hspace*{125mm}\rule{6pt}{6pt}

The goal of this paper is to provide reduced-order models that preserve
the positive realness and symmetry or block skew symmetry of the transfer matrix.
To apply RPRBT, the descriptor system of \rf{eqn:des1}
must be converted
to a state equation.
As the first step to obtain the state equation, consider the following
Weierstrass canonical form:
\begin{eqnarray}
\begin{aligned}
E_0 &= T_l
\left[
\begin{array}{cc}
I_q & 0 \\
0 & N
\end{array}
\right] T_r, \quad
A_0 = T_l
\left[
\begin{array}{cc}
J & 0 \\
0 & I_{n-q}
\end{array}
\right] T_r, \\
\end{aligned}
\label{eqn:weire}
\end{eqnarray}
where
$J$ is the Jordan form whose eigenvalues correspond to
the finite eigenvalues of the generalized eigenvalue problem $\left(E, A\right)$,
whereas
$N$ is a nilpotent whose eigenvalues are zero.
When $N^\mu = 0$, $N^{\mu'} \ne 0$, and $\mu' = \mu+1$,
$\mu$ is referred to as the index.
For RLC networks, the index is at most two \ct{phil}.

With the Weierstrass canonical form \rf{eqn:weire}, the relationships
$C_0 T_r^{-1} = [C_p~ C_{\infty}]$ and $T_l^{-1}B_0 = [B_p^T~ B_{\infty}^T]^T$
are defined.
Then, transfer matrix $G(s)$ is obtained as follows:
\begin{eqnarray}
\begin{aligned}
G(s) &= C_0 \left( s E_0 - A_0 \right)^{-1} B_0 \\
&=
\left[
\begin{array}{cc}
C_p & C_{\infty}
\end{array}
\right]
\left[
\begin{array}{cc}
\left(sI_q - J \right)^{-1} & \\
  & \left(sN - I_{n-q}\right)^{-1}
\end{array}
\right]
\left[
\begin{array}{c}
B_p \\
B_{\infty}
\end{array}
\right] \\
&= C_p \left( sI_q - J \right)^{-1} B_p + M_0 + s M_1,
\end{aligned}
\label{eqn:zs}
\end{eqnarray}
where $M_0 = - C_{\infty} B_{\infty}$ and $M_1 = - C_{\infty} N B_{\infty}$
are positive semidefinite \ct{phil}.
The terms $C_p \left( sI_q - J \right)^{-1} B_p + M_0$ in \rf{eqn:zs}
are the proper part of transfer matrix $G(s)$.
Note that $M_1$ becomes a zero matrix for an index-1 system.

%
%

\section{RPRBT for Index-1 Systems}

The descriptor system \rf{eqn:des1} is expressed in
a singular-value decomposition (SVD) canonical form
for conversion
to a state equation.
Although SVD is necessary, it is computationally
expensive for a large system.
Thus, rather than SVD, we use LDL factorization\footnote[1]{
Note that the MATLAB {\tt ldl}() function is available.}, which is
Cholesky-like factorization for semidefinite case:
$P^T \alpha P = L D L^T$,
where $\alpha$ is symmetric,
$P$ is a permutation, $L$ is a lower triangular matrix,
and $D$ is a diagonal matrix of rank less than full rank.
By applying LDL factorization to $A_{\mathcal C} \mathcal C A_{\mathcal C}$
and ${\mathcal L}$ in \rf{eqn:z}-\rf{eqn:h},
matrix $E_0$ can be expressed as follows:
\begin{eqnarray}
\begin{aligned}
E_0
&=
V
\left[
\begin{array}{ccc}
I_{r_1} & 0 & 0 \\
0 & -I_{r_2} & 0 \\
0 & 0 & 0
\end{array}
\right]
V^T
=
V
\left[
\begin{array}{cc}
I'_r & 0 \\
0 & 0
\end{array}
\right]
V^T,
\end{aligned}
\label{eqn:cong2}
\end{eqnarray}
where the ranks of $A_{\mathcal C} \mathcal C A_{\mathcal C}$
and $\mathcal L$ are assumed to be $r_1$ and $r_2$, respectively,
and $r = r_1 + r_2$.

Then, equation \rf{eqn:des1} is converted into the following:
\begin{eqnarray}
\begin{aligned}
\left[
\begin{array}{cc}
I'_r & 0 \\
0 & 0
\end{array}
\right]
\dfrac{d}{dt}
\left[
\begin{array}{c}
x_1(t) \\
x_2(t)
\end{array}
\right]
&=
\left[
\begin{array}{cc}
A_{11} & A_{12} \\
A_{21} & A_{22}
\end{array}
\right]
\left[
\begin{array}{c}
x_1(t) \\
x_2(t)
\end{array}
\right]
+
\left[
\begin{array}{c}
B_1 \\
B_2
\end{array}
\right] u(t), \\
z(t) &=
\left[
\begin{array}{cc}
C_1 & C_2
\end{array}
\right]
\left[
\begin{array}{c}
x_1(t) \\
x_2(t)
\end{array}
\right],
\end{aligned}
\label{eqn:svd}
\end{eqnarray}
which is considered an SVD canonical form of \rf{eqn:des1},
where $V^T x(t) = [x_1(t)^T ~ x_2(t)^T]^T$.

By assuming that $A_{22}$ is nonsingular, from \rf{eqn:svd},
we obtain the state equation as follows:
\begin{eqnarray}
\begin{aligned}
\frac{dx_1(t)}{dt} &= \hat{A}_0 x_1(t) + \hat{B}_0 u(t), \quad
z(t) = \hat{C}_0 x_1(t) + \hat{D}_0 u(t),
\end{aligned}
\label{eqn:des11}
\end{eqnarray}
where
\begin{eqnarray}
\begin{aligned}
\hat{A}_0 &= I'_r \left( A_{11} - A_{12} A_{22}^{-1} A_{21}\right) = I'_r A_s, \quad
\hat{B}_0 = I'_r \left( B_1 - A_{12} A_{22}^{-1} B_2 \right), \\
\hat{C}_0 &= C_1 - C_2 A_{22}^{-1} A_{21}, \quad
\hat{D}_0 = -C_2 A_{22}^{-1} B_2,
\end{aligned}
\label{eqn:sec2}
\end{eqnarray}
where $A_s$ is a symmetric matrix.
Here, the index is one if and only if
$A_{22}$ is nonsingular \ct{kata} and the descriptor system
is always converted to a state equation.

The following two AREs are solved for PRBT:
\begin{eqnarray}
A^T X + X A + XBB^TX + C^TC = 0,
\label{eqn:ric1} \\
A Y + Y A^T + YC^TCY + BB^T = 0,
\label{eqn:ric2}
\end{eqnarray}
where $DD^T = ( \hat{D}_0 + \hat{D}_0^T )^{-1}$,
$B = \hat{B}_0 D$,
$C = D^T \hat{C}_0$, and $A = \hat{A}_0 - BC$.

Prior to applying PRBT,
we prove that the solution of ARE
\rf{eqn:ric2} can be expressed by
that of \rf{eqn:ric1}.
\begin{theorem}
\label{th:are}
The solution of ARE
\rf{eqn:ric2} is expressed as $Y = I'_r X I'_r$, with that of \rf{eqn:ric1}.
\end{theorem} \\
{\it Proof.}~See Appendix A. \\
\hspace*{125mm}\rule{6pt}{6pt}

Here, consider the similarity transformations as
$\tilde{A}_0 = T^{-1} \hat{A}_0 T$,
$\tilde{B}_0 = T^{-1} \hat{B}_0$, and
$\tilde{C}_0 = \hat{C}_0 T$.
The solutions of the transformed AREs are
both diagonalized as follows:
\begin{eqnarray}
T^T X T = \Sigma = T^{-1} Y T^{-T}.
\label{eqn:gram}
\end{eqnarray}
The diagonal elements of $\Sigma$ are called
Hankel singular values.
Note that small values exhibit weak effects on the input--output
behavior \ct{phil}; therefore, components corresponding
to small Hankel singular values can be removed.

From theorem \ref{th:are}, $Y = I'_r Q Q^T I'_r$ is obtained using
Cholesky factorization $X=Q Q^T$.
Moreover, using eigenvalue decomposition, we obtain
$Q^T I'_r Q = U \Sigma U^T = U \left|\Sigma\right| S U^T$,
where $S$ is a diagonal matrix whose diagonal elements
are the signs of $\Sigma$. Then,
the transformation matrices are obtained by $T = I_r' Q U S \left|\Sigma\right|^{-1/2}$
and $T^{-1} = \left|\Sigma\right|^{-1/2} U^T Q^T$.
Assuming the absolute values of the diagonal elements of $\Sigma$
are arranged in descending order,
matrices $\Sigma$, $S$, and $U$
are partitioned as follows:
\begin{eqnarray}
\begin{aligned}
\Sigma &=
\left[
\begin{array}{cc}
\Sigma_1 & 0 \\
0 & \Sigma_2
\end{array}
\right], \quad
S =
\left[
\begin{array}{cc}
S_1 & 0 \\
0 & S_2
\end{array}
\right], \quad
U &=
\left[
\begin{array}{cc}
U_1 & U_2
\end{array}
\right].
\end{aligned}
\label{eqn:bt_rlc11}
\end{eqnarray}
With the Hankel singular values included
in only $\Sigma_1$,
the balanced realization is obtained by
\begin{eqnarray}
\begin{aligned}
\tilde{A}_{11} &=
\left|\Sigma_1\right|^{-1/2} U_1^T Q^T
I'_r A_s I'_r Q U_1 \left|\Sigma_1\right|^{-1/2} S_1 = \tilde{A}_s S_1, \\
\tilde{B}_{1} &= \left|\Sigma_1\right|^{-1/2} U_1^T Q^T \hat{B}_0, \quad
\tilde{C}_{1} = \hat{C}_0 I'_r Q U_1
S_1 \left|\Sigma_1\right|^{-1/2},
\end{aligned}
\label{eqn:bt_rlc12}
\end{eqnarray}
where $\tilde{A}_s$ is a symmetric matrix.

Then, the transfer matrix of the reduced-order model
is represented as:
\begin{eqnarray}
\tilde{G}(s)
&=& \tilde{C}_1 S_1 \left( s S_1 - \tilde{A}_s \right)^{-1} \tilde{B}_1
+ \hat{D}_0.
\label{eqn:rs_rlc2}
\end{eqnarray}
Thus, we obtain the following theorem.

\begin{theorem}
\label{th:res1}
The reciprocity of the descriptor system \rf{eqn:des1} is
preserved after applying PRBT.
\end{theorem} \\
{\it Proof.}~The reduced-order impedance, admittance,
and hybrid matrices are expressed as follows:
\begin{eqnarray}
\begin{aligned}
Z(s) &= \tilde{B}_1^T \left( s S_1 - \tilde{A}_s \right)^{-1} \tilde{B}_1
- B_2^T A_{22}^{-1} B_2, \\
Y(s) &= -\tilde{B}_1^T \left( s S_1 - \tilde{A}_s \right)^{-1} \tilde{B}_1
+ B_2^T A_{22}^{-1} B_2, \\
H(s) &= I^\circ_m \tilde{B}_1^T \left( s S_1 - \tilde{A}_s \right)^{-1} \tilde{B}_1
- I^\circ_m B_2 A_{22}^{-1} B_2,
\end{aligned}
\end{eqnarray}
where $Z(s)$, $Y(s)$, and $H(s)$ are the impedance, admittance,
and hybrid matrices, respectively.
As matrix $( s S_1 - \tilde{A}_s )^{-1}$ is symmetric,
$Z(s)$ and $Y(s)$ are symmetric and
$H(s)$ is block skew symmetric. \\
\hspace*{125mm}\rule{6pt}{6pt}

In this paper, the reciprocity-preserving PRBT algorithm
for an index-1 system is called {\bf RPRBT-1}.
\begin{algorithm}
\caption{~\bf{RPRBT-1}}
\begin{enumerate}
\setlength{\parskip}{0cm}
\item Solve \rf{eqn:ric1} for $X$.
\item Compute the Cholesky factor as $X = Q Q^T$.
\item Compute the eigenvalue decomposition
as $U \left|\Sigma\right| S U^T = Q^T I'_r Q$.
\item Compute the transformation matrices as
$T = I_r' Q U S \left|\Sigma\right|^{-1/2}$
and $T^{-1} = \left|\Sigma\right|^{-1/2} U^T Q^T$.
\item Compute the reduced-order matrices as
$\tilde{A}_0 = T^{-1} \hat{A}_0 T$,
$\tilde{B}_0 = T^{-1} \hat{B}_0$, and
$\tilde{C}_0 = \hat{C}_0 T$.
\item Partition $\tilde{A}_0$, $\tilde{B}_0$,
and $\tilde{C}_0$
as
\begin{eqnarray*}
\tilde{A}_0 &=&
\left[
\begin{array}{cc}
\tilde{A}_{11} & \tilde{A}_{12} \\
\tilde{A}_{21} & \tilde{A}_{22}
\end{array}
\right], \quad
\tilde{B}_0 =
\left[
\begin{array}{c}
\tilde{B}_{1} \\
\tilde{B}_{2}
\end{array}
\right], \quad
\tilde{C}_0 =
\left[
\begin{array}{cc}
\tilde{C}_{1} & \tilde{C}_{2}
\end{array}
\right].
\end{eqnarray*}
\item Truncate $\tilde{A}_0$, $\tilde{B}_0$, and $\tilde{C}_0$
to form the reduced realization $( \tilde{A}_{11},
\tilde{B}_1, \tilde{C}_1, \hat{D}_0 )$;
the reduced-order transfer matrix is then obtained by
\rf{eqn:rs_rlc2}.
\end{enumerate}
\end{algorithm}
\\
Note that Cholesky factorization
of step 2 is not necessary because the Cholesky
factor is obtained directly by ARE solvers (e.g., \ct{benner2}, \ct{reis}).

In the SVD canonical form, we cannot guarantee that matrix $\hat{A}_0$ in \rf{eqn:sec2}
is nonsingular; thus,
there does not exist a value of transfer matrix
$\hat{C}_0 ( s I_r - \hat{A}_0 )^{-1} \hat{B}_0 + \hat{D}_0$ at $s=0$,
even if the original
$C_0 \left( s E_0 - A_0 \right)^{-1} B_0 + D_0$ has a value at $s=0$.
This occurs due to the introduction of the SVD canonical form.
Therefore, we must eliminate this artifact.
With a permutation matrix $P$, the strict proper part of the
transfer matrix is rewritten as follows:
\begin{eqnarray}
\begin{aligned}
G_p(s) &= \hat{C}_0 \left( s I_r - I'_r A_s \right)^{-1} \hat{B}_0 \\
  &= \hat{C}_0 P^T
\left( s I''_r -
\left[
\begin{array}{cc}
\bar{A}_{s,1} & 0 \\
0 & 0
\end{array}
\right]
\right)^{-1} P I'_r \hat{B}_0 \\
&=
\left[
\begin{array}{cc}
\hat{C}_1 & \hat{C}_2
\end{array}
\right]
\left[
\begin{array}{cc}
s I''_1 - \hat{A}_{s,1} & 0 \\
0 & s I''_2
\end{array}
\right]^{-1}
\left[
\begin{array}{c}
\hat{B}_1 \\
\hat{B}_2
\end{array}
\right] \\
&= \hat{C}_1 \left( s I''_1 - \hat{A}_{s,1} \right)^{-1} \hat{B}_1
+ \hat{C}_2 \left( s I''_2 \right)^{-1} \hat{B}_2.
\end{aligned}
\label{eqn:elm}
\end{eqnarray}
If there exists a value of the original system at $s=0$,
the second term of \rf{eqn:elm} must be eliminated;
thus, RPRBT is applied to
the realization $( I''_1 \hat{A}_{s,1}, I''_1 \hat{B}_1,
\hat{C}_1, \hat{D}_0 )$.

%
%

\section{RPRBT for Index-2 Systems}

\subsection{GARE}

When matrix $A_{22}$ in \rf{eqn:svd} is singular,
the index becomes two for passive RLC networks.
Therefore, all terms of \rf{eqn:zs} must be calculated.
Then, we introduce the right and left spectral projectors
associated with the deflated invariant subspace of
matrix pencil $\lambda E_0 - A_0$.
To obtain the spectral projectors,
LDL decomposition is applied to $A_{22}$ and
the following relationship is obtained:
\begin{eqnarray}
A_{22} &=&
F
\left[
\begin{array}{cc}
S_{\bar{r}} & 0 \\
0 & 0
\end{array}
\right]
F^T,
\label{eqn:ldl2}
\end{eqnarray}
where matrix $S_{\bar{r}}$ is nonsingular.
By expressing $F^T x_2(t) = y(t) = [y_1(t)^T~ y_2(t)^T]^T$,
we rewrite \rf{eqn:svd} as follows:
\begin{eqnarray}
\begin{aligned}
I'_r \frac{d x_1(t)}{dt} &= A_{11} x_1(t) + A_{12,1} y_1(t)
+ A_{12,2} y_2(t) + B_1 u(t), \\
0 &=
\left[
\begin{array}{c}
A_{21,1} \\
A_{21,2}
\end{array}
\right]
+
\left[
\begin{array}{c}
S_{\bar{r}} y_1(t) \\
0
\end{array}
\right]
+
\left[
\begin{array}{c}
B_{2,1} \\
B_{2,2}
\end{array}
\right]
u(t), \\
z(t) &= C_1 x_1(t) + C_{2,1} y_1(t) + C_{2,2} y_2(t).
\end{aligned}
\label{eqn:svd3}
\end{eqnarray}
By eliminating $y_1(t)$ in \rf{eqn:svd3}, the
following relationship is obtained:
\begin{eqnarray}
\begin{aligned}
\left[
\begin{array}{cc}
I'_r & 0 \\
0 & 0
\end{array}
\right]
\frac{d}{dt}
\left[
\begin{array}{c}
x_1(t) \\
y_2(t)
\end{array}
\right]
&=
\left[
\begin{array}{cc}
\bar{A}_{11} & \bar{A}_{12} \\
\bar{A}_{21} & 0
\end{array}
\right]
\left[
\begin{array}{c}
x_1(t) \\
y_2(t)
\end{array}
\right]
+
\left[
\begin{array}{c}
\bar{B}_1 \\
\bar{B}_2
\end{array}
\right] u(t), \\
z(t) &= \bar{C}_1 x_1(t) + \bar{C}_2 y_2(t) + \bar{D}_0 u(t),
\end{aligned}
\label{eqn:stokes}
\end{eqnarray}
where $\bar{A}_{11} = A_{11} - A_{12,1} S_{\bar{r}}^{-1} A_{21,1}$,
$\bar{A}_{12} = A_{12,2}$,
$\bar{A}_{21} = A_{21,2}$,
$\bar{B}_1 = B_1 - A_{12,1} S_{\bar{r}}^{-1} B_{2,1}$,
$\bar{B}_2 = B_{2,2}$,
$\bar{C}_1 = C_1 - C_{2,1} S_{\bar{r}}^{-1} A_{21,1}$,
$\bar{C}_2 = C_{2,2}$, and
$\bar{D}_0 = - C_{2,1} S_{\bar{r}}^{-1} B_{2,1}$.
Equation \rf{eqn:stokes} is referred to as a Stokes-type index-2 system
whose left and right spectral projectors $P_l$ and $P_r$
are explicitly written as follows:
\begin{eqnarray}
\begin{aligned}
P_l &=
\left[
\begin{array}{cc}
\Pi_l & -\Pi_l \bar{A}_{11} I'_r \bar{A}_{12} \left( \bar{A}_{21} I'_r \bar{A}_{12} \right)^{-1} \\
0 & 0
\end{array}
\right], \\
P_r &=
\left[
\begin{array}{cc}
\Pi_r & 0 \\
-\left( \bar{A}_{21} I'_r \bar{A}_{12} \right)^{-1} \bar{A}_{21} I'_r \bar{A}_{11} \Pi_r & 0
\end{array}
\right],
\end{aligned}
\label{eqn:proj}
\end{eqnarray}
where
$\Pi_l = I_r - \bar{A}_{12} \left( \bar{A}_{21} I'_r \bar{A}_{12} \right)^{-1}
\bar{A}_{21} I'_r$
and $\Pi_r = I_r - I'_r \bar{A}_{12} \left( \bar{A}_{21} I'_r
\bar{A}_{12} \right)^{-1} \bar{A}_{21}$.
$\Pi_l$ is a projector onto the kernel of $\bar{A}_{21} I'_r$
along the image of $\bar{A}_{12}$.

To apply PRBT, the descriptor system \rf{eqn:stokes} is rewritten as:
\begin{eqnarray}
\begin{aligned}
\bar{E}_0 \dfrac{d \bar{x}(t)}{dt} &= \bar{A}_0 \bar{x}(t)
+ \bar{B}_0 u(t), \quad
z(t) = \bar{C}_0 \bar{x}(t) + \bar{D}_0 u(t).
\end{aligned}
\label{eqn:des2}
\end{eqnarray}
The right and left projectors satisfy the following relationships:
\begin{eqnarray}
P_r = T_r^{-1}
\left[
\begin{array}{cc}
I_r & 0 \\
0 & 0
\end{array}
\right]
T_r, \quad
P_l = T_l
\left[
\begin{array}{cc}
I_r & 0 \\
0 & 0
\end{array}
\right]
T_l^{-1}.
\label{eqn:sp2}
\end{eqnarray}
Thus, we can express $M_0$ and $M_1$ of \rf{eqn:zs} as follows:
\begin{eqnarray}
\begin{aligned}
M_0 &= - \bar{C}_0 \left( I_{n-\bar{r}} - P_r \right) \bar{A}_0^{-1}
\left( I_{n-\bar{r}} - P_l \right) \bar{B}_0 + \bar{D}_0, \\
M_1 &= - \bar{C}_0 \bar{A}_0^{-1} \left( I_{n-\bar{r}} - P_l \right)
\bar{E}_0 \left( I_{n-\bar{r}} - P_r \right) \bar{A}_0^{-1} \bar{B}_0,
\end{aligned}
\label{eqn:m0m1}
\end{eqnarray}
which are proven in Appendix B.
If $\bar{E}_0$ and $\bar{A}_0$ are symmetric, $P_l = P_r^T$ \ct{wong3}.
As the two transforms \rf{eqn:cong2} and \rf{eqn:ldl2} do not
break the symmetry of the original descriptor system,
$P_l = P_r^T$ holds for \rf{eqn:proj}.
Therefore, $M_0$ and $M_1$ in \rf{eqn:m0m1} are symmetric to
the impedance and admittance matrices, and are
block skew symmetric to the hybrid matrix.

In PRBT, the following dual GAREs are solved \ct{stykel4}.
\begin{eqnarray}
&&A_1^T X E + E^T X A_1 + E^T X BB^T X E + P_r^T C^T C P_r = 0,
\label{eqn:gric1} \\
&&A_2 Y E^T + E Y A_2^T + E Y C^TC Y E^T + P_l B B^T P_l^T = 0,
\label{eqn:gric2}
\end{eqnarray}
where $E=\bar{E}_0$, $DD^T = \left( M_0 + M_0^T \right)^{-1}$,
$A_1 = \bar{A}_0 - \bar{B}_0 D D^T \bar{C}_0 P_r$,
$A_2 = \bar{A}_0 - P_l \bar{B}_0 D D^T \bar{C}_0$,
$C = D^T \bar{C}_0$, and $B = \bar{B}_0 D$.
Then, the following theorem holds for
GAREs \rf{eqn:gric1} and \rf{eqn:gric2}.
\begin{theorem}
\label{th:gare}
The solutions of the dual GAREs \rf{eqn:gric1} and \rf{eqn:gric2} are
equal, i.e., $Y = X$.
\end{theorem} \\
{\it Proof:}~See Appendix C. \\
\hspace*{125mm}\rule{6pt}{6pt}

From theorem \ref{th:gare}, $X = Q Q^T$ is obtained using the Cholesky factorization
$Y = Q Q^T$.
In addition, using eigenvalue decomposition, we obtain
$Q^T \bar{E}_0 Q = U \left|\Sigma\right| S U^T$,
where $S$ is a diagonal matrix, the diagonal elements of which
are the signs of $\Sigma$.
Then, the transformation matrices are obtained by
$\bar{T} = Q U S \left|\Sigma\right|^{-1/2}$
and $\bar{W} = \left|\Sigma\right|^{-1/2} U^T Q^T$.
Assuming the absolute values of the diagonal elements of $\Sigma$
are arranged in descending order,
matrices $\Sigma$, $S$, and $U$
are partitioned as follows:
\begin{eqnarray}
\begin{aligned}
\Sigma &=
\left[
\begin{array}{cc}
\Sigma_1 & 0 \\
0 & \Sigma_2
\end{array}
\right], \quad
S =
\left[
\begin{array}{cc}
S_1 & 0 \\
0 & S_2
\end{array}
\right],
\quad
U &=
\left[
\begin{array}{cc}
U_1 & U_2
\end{array}
\right].
\end{aligned}
\label{eqn:bt_rlc21}
\end{eqnarray}
With the Hankel singular values included
in only $\Sigma_1$,
the balanced realization is obtained by
\begin{eqnarray}
\begin{aligned}
\tilde{E}_{11} &=
\left|\Sigma_1\right|^{-1/2} U_1^T Q^T \bar{E}_0 Q U_1 S_1 \left|\Sigma_1\right|^{-1/2}
= \tilde{E}_s S_1, \\
\tilde{A}_{11} &=
\left|\Sigma_1\right|^{-1/2} U_1^T Q^T \bar{A}_0 Q U_1 S_1 \left|\Sigma_1\right|^{-1/2}
= \tilde{A}_s S_1, \\
\tilde{B}_{1} &= \left|\Sigma_1\right|^{-1/2} U_1^T Q^T \bar{B}_0, \quad
\tilde{C}_{1} = \bar{C}_0 Q U_1 S_1 \left|\Sigma_1\right|^{-1/2},
\end{aligned}
\label{eqn:bt_rlc22}
\end{eqnarray}
where $\tilde{E}_s$ and $\tilde{A}_s$ are symmetric.
The transfer matrix of the reduced-order model
is expressed as follows:
\begin{eqnarray}
\tilde{G}(s)
&=& \tilde{C}_1 S_1 \left( s \tilde{E}_s
- \tilde{A}_s \right)^{-1} \tilde{B}_1
+ M_0 + s M_1.
\label{eqn:rs_rlc3}
\end{eqnarray}
The reciprocity-preserving PRBT algorithm
for an index-2 system is called {\bf RPRBT-2}.
\begin{algorithm}
\caption{~\bf{RPRBT-2}}
\begin{enumerate}
\setlength{\parskip}{0cm}
\item Solve \rf{eqn:gric1} for $X$.
\item Compute the Cholesky factor as $X = Q Q^T$.
\item Apply eigenvalue decomposition
as $U \left|\Sigma\right| S U^T = Q^T E_0 Q$.
\item Compute the transformation matrices as
$\bar{T} = Q U S \left|\Sigma\right|^{-1/2}$ and
$\bar{W} = \left|\Sigma\right|^{-1/2} U^T Q^T$.
\item Compute the reduced-order matrices as
$\tilde{E}_0 = \bar{W} \bar{E}_0 \bar{T}$,
$\tilde{A}_0 = \bar{W} \bar{A}_0 \bar{T}$,
$\tilde{B}_0 = \bar{W} \bar{B}_0$, and
$\tilde{C}_0 = \bar{C}_0 \bar{T}$.
\item Partition $\tilde{E}_0$, $\tilde{A}_0$,
$\tilde{B}_0$, and $\tilde{C}_0$
as
\begin{eqnarray*}
\tilde{E}_0 &=&
\left[
\begin{array}{cc}
\tilde{E}_{11} & \tilde{E}_{12} \\
\tilde{E}_{21} & \tilde{E}_{22}
\end{array}
\right], \quad
\tilde{A}_0 =
\left[
\begin{array}{cc}
\tilde{A}_{11} & \tilde{A}_{12} \\
\tilde{A}_{21} & \tilde{A}_{22}
\end{array}
\right], \quad
\tilde{B}_0 =
\left[
\begin{array}{c}
\tilde{B}_{1} \\
\tilde{B}_{2}
\end{array}
\right], \\
\tilde{C}_0 &=&
\left[
\begin{array}{cc}
\tilde{C}_{1} & \tilde{C}_{2}
\end{array}
\right].
\end{eqnarray*}
\item Truncate $\tilde{E}_0$, $\tilde{A}_0$, $\tilde{B}_0$, and
$\tilde{C}_0$
to form the reduced realization
$( \tilde{E}_{11}, \tilde{A}_{11},
\tilde{B}_1, \tilde{C}_1, M_0, M_1 )$;
the reduced-order transfer function is obtained by
\rf{eqn:rs_rlc3}.
\end{enumerate}
\end{algorithm}
Note that
Cholesky factorization of step 2 is not required for {\bf RPRBT-1}
because it is obtained by a GARE solver.

Then, we obtain the following theorem.

\begin{theorem}
The reciprocity of the descriptor system with \rf{eqn:des1} is
preserved after applying {\bf RPRBT-2}, even if the systems are index-2.
\end{theorem} \\
{\it Proof.}~%
From the symmetry of $( s \tilde{E}_s - \tilde{A}_s )^{-1}$,
$Z(s)$ and $Y(s)$ are symmetric, and $H(s)$
is block skew symmetric. \\
\hspace*{125mm}\rule{6pt}{6pt}

For an index-1 system, such as $S_r \in R^{n-r,n-r}$ in \rf{eqn:ldl2},
equation \rf{eqn:stokes} is expressed as follows:
\begin{eqnarray}
\begin{aligned}
I'_r \frac{d x_1(t)}{dt}
&=
\bar{A}_{11} x_1(t) +
\bar{B}_1 u(t), \\
z(t) &= \bar{C}_1 x_1(t) + \bar{D}_0 u(t).
\end{aligned}
\label{eqn:stokes2}
\end{eqnarray}
By inputting
$\bar{E}_0 = I'_r$,
$\bar{A}_0 = \bar{A}_{11}$,
$\bar{B}_0 = \bar{B}_1$,
$\bar{C}_0 = \bar{C}_1$,
$P_r = P_l = I_r$, and $M_0 = \bar{D}_0$ in \rf{eqn:des2},
{\bf SPRBT-2} can be applied to the index-1 system.


\subsection{ARE}

In the previous subsection, the reduced-order model
was obtained via the GARE.
As AREs have been studied more than GAREs,
it is preferable to define an ARE for the proper part
of \rf{eqn:zs} and apply PRBT to this part.
We define the ARE beginning from \rf{eqn:stokes} and applying the projector-based
methods \ct{proj1} and \ct{proj2}.

From the second block of the first equation of \rf{eqn:stokes},
there exists a special solution:
\begin{eqnarray}
x_{1g}(t) &=& - I'_r \bar{A}_{12}
\left( \bar{A}_{21} I_r' \bar{A}_{12} \right)^{-1} \bar{B}_2 u(t).
\label{eqn:ss}
\end{eqnarray}
Representing $x_1(t) = x_{10}(t) + x_{1g}(t)$ and using the second equation
of \rf{eqn:stokes}, we obtain the following:
\begin{eqnarray}
\bar{A}_{21} x_{10}(t) &=& 0.
\label{eqn:a21}
\end{eqnarray}
From the first equation of \rf{eqn:stokes}, we obtain:
\begin{eqnarray}
\dot{x}_{10}(t) &=& I'_r \bar{A}_{11} x_{10}(t) + I'_r \bar{A}_{12} y_{2}(t)
+ I'_r \bar{A}_{11} x_{1g}(t) + I'_r \bar{B}_{1} u(t) - \dot{x}_{1g}(t),
\label{eqn:x01}
\end{eqnarray}
where $\dot{x}$ indicates $dx/dt$.
Using \rf{eqn:a21} and \rf{eqn:x01}, $y_2(t)$ is expressed by:
\begin{eqnarray}
y_2(t) &=& - \left(\bar{A}_{21} I'_r \bar{A}_{12} \right)^{-1}
\left\{ \bar{A}_{21} I'_r \bar{A}_{11}
\left( x_{10}(t) + x_{1g}(t) \right) \right.
\nonumber \\
&& \left. + \bar{A}_{21} I'_r \bar{B}_{1} u(t)
- \bar{A}_{21} \dot{x}_{1g}(t) \right\}.
\label{eqn:y2}
\end{eqnarray}
By eliminating $y_2(t)$ in \rf{eqn:x01} and using the left projector
$\Pi_l$ of \rf{eqn:proj}, we obtain the following:
\begin{eqnarray}
\dot{x}_{10}(t) &=& I'_r \Pi_l A_{11} x_{10}(t) +
I'_r \Pi_l \left\{ \bar{B}_1 - \bar{A}_{11} I'_r \bar{A}_{12}
\left( \bar{A}_{21} I'_r \bar{A}_{12} \right)^{-1}
\bar{B}_2 \right\} u(t),
\label{eqn:x01_2}
\end{eqnarray}
where $- \dot{x}_{1g}(t) +
I'_r \bar{A}_{12} \left( \bar{A}_{21} I'_r \bar{A}_{12} \right)^{-1}
\bar{A}_{21} \dot{x}_{1g}(t) = 0$ is used.

With \rf{eqn:a21} and the right projector $\Pi_r$,
$\Pi_r x_{10}(t) = x_{10}(t)$ holds.
Therefore, the Laplace transform of \rf{eqn:x01_2} is expressed as follows:
\begin{eqnarray}
X_{01}(s) &=& \left(s I_r - I'_r \Pi_l A_{11} \Pi_r \right)^{-1} \nonumber \\
&&
\hspace*{-10mm}
\times
I'_r \Pi_l \left\{ \bar{B}_1 - \bar{A}_{11} I'_r \bar{A}_{12}
\left( \bar{A}_{21} I'_r \bar{A}_{12} \right)^{-1}
\bar{B}_2 \right\} U(s),
\label{eqn:x01_3}
\end{eqnarray}
where
$X_{01}(s)$ and $U(s)$ are the Laplace transforms of $x_{01}(t)$ and $u(t)$,
respectively.
Then,
transforming the second equation of \rf{eqn:stokes}
into the Laplace domain,
transfer function $G(s)$ is expressed as:
\begin{eqnarray}
G(s) &=& C_p \left( s I_r - J_p \right)^{-1} B_p + M_0 + s M_1,
\label{eqn:ge}
\end{eqnarray}
where
\begin{eqnarray}
\begin{aligned}
J_p &= I'_r \Pi_l \bar{A}_{11} \Pi_r, \\
C_p &= \left\{ \bar{C}_1 -
\bar{C}_2 \left( \bar{A}_{21} I'_r \bar{A}_{12} \right)^{-1}
\bar{A}_{21} I'_r \bar{A}_{11} \right\} \Pi_r, \\
B_p &= I'_r \Pi_l \left\{ \bar{B}_1 - \bar{A}_{11} I'_r \bar{A}_{12}
\left( \bar{A}_{21} I'_r \bar{A}_{12} \right)^{-1} \bar{B}_2 \right\}, \\
M_0 &= - \bar{C}_1 I'_r \bar{A}_{12}
\left( \bar{A}_{21} I'_r \bar{A}_{12} \right)^{-1} \bar{B}_2
- \bar{C}_2 \left( \bar{A}_{21} I'_r \bar{A}_{12} \right)^{-1}
\bar{A}_{21} I'_r \bar{B}_1 \\
&
\hspace*{10mm}
+ \bar{C}_2 \left( \bar{A}_{21} I'_r \bar{A}_{12} \right)^{-1}
\bar{A}_{21} I'_r \bar{A}_{11} I'_r \bar{A}_{12}
\left( \bar{A}_{21} I'_r \bar{A}_{12} \right)^{-1} \bar{B}_2 + \bar{D}_0, \\
M_1 &= - \bar{C}_2 \left( \bar{A}_{21} I'_r \bar{A}_{12} \right)^{-1} \bar{B}_2.
\end{aligned}
\label{eqn:coef}
\end{eqnarray}
The coefficient matrices of \rf{eqn:coef} are similar to
\rf{eqn:sec2}; thus, {\bf RPRBT-1} is applied to the proper part,
and the reduced-order model is obtained as follows:
\begin{eqnarray}
\tilde{G}(s)
&=& \tilde{C}_1 S_1 \left( s S_1 - \tilde{A}_s \right)^{-1} \tilde{B}_1
+ M_0 + s M_1.
\label{eqn:rs_rlc4}
\end{eqnarray}
Then, the reciprocity of the reduced-order model is preserved,
and the following theorem is obtained without proof.
\begin{theorem}
\label{th:res2-1}
The reciprocity of the descriptor system with \rf{eqn:des1} is
preserved after applying {\bf RPRBT-1}
to the proper part of \rf{eqn:ge}, even if the system is
index-2.
\end{theorem}
%

%
%

\section{RADI}

The RADI used to solve GARE \rf{eqn:gare} is described in Algorithm 3.
The RADI used to solve AREs is provided with $E =I$ in Algorithm 3;
thus, it is a special case of Algorithm 3.
The approximate feedback matrix $K = E^* X B$ is introduced
to apply the Sherman-Morrison-Woodbury (SMW) formula, which accelerates
the algorithm.
By defining the solution at the $k$-th while loop as $X_k$
and matrices $V$ and $\tilde{Y}$ as $V_k$ and $\tilde{Y}_k$, respectively,
the GARE solution is expressed as follows:
\begin{eqnarray}
X_k = \sum_{i=1}^k V_i \tilde{Y}_i^{-1} V_i^*,
\label{eqn:sol_radi}
\end{eqnarray}
where $X_0 = 0$.
Moreover, the solution is expressed as $X_k = Z_k Y_k^{-1} Z_k^*$,
where $Z_k$ and $Y_k$ are matrices $Z$ and $Y$ at the $k$-th loop, respectively.

In step 2 of {\bf RPRBT-1} and {\bf RPRBT-2}, we calculate the Cholesky factor $Q$,
which is obtained as follows:
\begin{eqnarray}
Q = \left[ \mbox{Re} (Z_k Y_{k,h}^{-1})~\mbox{Im}(Z_k Y_{k,h}^{-1}) \right],
\label{eqn:cf_radi}
\end{eqnarray}
where $Y_k^{-1} = Y_{k,h}^{-1} Y_{k,h}^{-*}$. As $Y$ is a block diagonal matrix,
matrix $Y_{k,h}^{-1}$ can be calculated efficiently.
Algorithm 4 is a fundamental implementation of RADI.
A more efficient implementation that avoids complex arithmetic
for complex conjugate shifts is provided in the literature \ct{benner2}.

\begin{algorithm}
\caption{~RADI}
\begin{algorithmic}
\STATE $R = C^T, K = 0, Y = [~]$;
\WHILE{$\| R^* R \| \geq {\rm tol} \cdot \| C C^T \|$}
\STATE Obtain the next shift $\sigma$;
\IF{first pass through the loop}
\STATE $Z = V = \sqrt{-2 {\rm Re}(\sigma)} (A^T + \sigma E^T)^{-1} R$;
\ELSE
\STATE $V = \sqrt{-2 {\rm Re}(\sigma)} (A^T + K B^T + \sigma E^T)^{-1} R$; // Use SMW
\STATE $Z = [Z~V];$
\ENDIF
\STATE $\tilde{Y} = I + \frac{1}{2{\rm Re}(\sigma)}(V^* B) (V^* B)^*;
Y = \left[
\begin{array}{cc}
Y & \\
  & \tilde{Y}
\end{array} \right] $;
\STATE $R = R + \sqrt{-2 {\rm Re}(\sigma)} E^T V \tilde{Y}^{-1}$;
\STATE $K = K + E^T V \tilde{Y}^{-1} V^* B$;
\ENDWHILE
\end{algorithmic}
\end{algorithm}

The performance of ADI depends on shifts that have negative real parts.
A thorough analysis of shift selection is provided in the literature \ct{benner}.
Note that shift selections are effective for obtaining a better
ARE solution, i.e.,
a low-rank solution has small ARE residual error.
Real and complex conjugate eigenvalues of a Hamiltonian matrix
are used as shifts
to reduce the ARE residual error effectively.
However, shift selections are insufficient for MOR of
electrical circuits \ct{tanji}.
The goal of MOR for electrical circuits is to obtain
an accurate low-frequency model.
When shifts are selected such that the ARE residual error
is reduced considerably, eigenvalues with small radius from the origin
(expressed as small eigenvalues throughout this paper)
tend not to be selected.
However, the small eigenvalues contribute to
model accuracy at low frequencies.

To solve a large-scale ARE,
eigenvalues are approximated using a Krylov subspace method,
such as the Arnoldi method.
Eigenvalues with large radius from the origin (expressed as large eigenvalues)
are obtained easily by the Krylov subspace method,
and small eigenvalues are obtained by applying
the Krylov subspace method to the inverse Hamiltonian matrix.
However, the shift selections provided in the literature \ct{tanji}
could not find suitable small eigenvalues.
Thus, a small negative constant value is used to compensate
model accuracy at low frequencies.

To improve the shift selections,
we first calculate the eigenvalues of the inverse Hamiltonian matrix
using the Krylov subspace method. Next, with all the ones with negative real
parts used to solve the ARE,
the reciprocal values are used as the shifts of RADI to solve the ARE.
On the other hand,
the generalized eigenvalue problem associated with GARE is expressed as follows:
$\lambda {\cal E}x = {\cal H}x$, where
\begin{eqnarray*}
{\cal H} &=&
\left[
\begin{array}{cc}
A & B B^T \\
-C^T C & -A^T
\end{array}
\right], \quad
{\cal E} = {\rm diag}\{E, E^*\}.
\end{eqnarray*}
This equation is rewritten as ${\cal H}^{-1} {\cal E} x = (1/\lambda) x$.
Then, the eigenvalues of matrix ${\cal H}^{-1} {\cal E}$
are calculated using the Krylov subspace method, and the reciprocal values
are used as the shifts of RADI to solve the GARE.
The shift computation for RADI is described in Algorithm 4.
\begin{algorithm}
\caption{~\bf{Shift Computation for RADI}}
\begin{enumerate}
\item Calculate eigenvalues of inverse Hamiltonian matrix ${\cal H}^{-1}$
to solve the ARE and matrix ${\cal H}^{-1} {\cal E}$ to solve the GARE
using the Krylov subspace method.
\item Select eigenvalues with negative real parts.
\item Calculate the reciprocal values.
\end{enumerate}
\end{algorithm}

The eigenvalues obtained by Krylov subspace methods
are extreme eigenvalues of a matrix.
In other words, we obtain both small and large values.
The large values are effective for reducing ARE or GARE residual error,
which improves model accuracy at high frequencies.

%
%
\section{Results}


\subsection{ARE Solution for Index-1 System}

Consider the descriptor system \rf{eqn:des1} with the following
coefficient matrices:
\begin{eqnarray*}
\begin{aligned}
E_0 &=
\left[
\begin{array}{rrrrrrr}
 1 & 0 & 0 & 0 & 0 & 0 & 0 \\
 0 & 0 & 0 & 0 & 0 & 0 & 0 \\
 0 & 0 & 1 & 0 & 0 & 0 & 0 \\
 0 & 0 & 0 & 0 & 0 & 0 & 0 \\
 0 & 0 & 0 & 0 & 0 & 0 & 0 \\
 0 & 0 & 0 & 0 & 0 & 1 & 0 \\
 0 & 0 & 0 & 0 & 0 & 0 & 1
\end{array}
\right], \quad
-A_0 =
\left[
\begin{array}{rrrrrrr}
  2 & -1 & 0 & 0 & 0 & 0 & 0 \\
 -1 & 1 & 0 & 0 & 0 & 1 & 0 \\
  0 & 0 & 1 & -1 & 0 & -1 & 0 \\
  0 & 0 & -1 & 1 & 0 & 0 & 1 \\
  0 & 0 & 0 & 0 & 0 & 0 & -1 \\
  0 & 1 & -1 & 0 & 0 & 0 & 0 \\
  0 & 0 & 0 & 1 & -1 & 0 & 0
\end{array}
\right], \quad \\
C_0 &=
\left[
\begin{array}{rrrrrrr}
-1 & 0 & 0 & 0 & 0 & 0 & 0 \\
0 & 0 & 0 & 0 & -1 & 0 & 0
\end{array}
\right], \quad
B_0 = C_0^T, \quad D_0 = 0.
\end{aligned}
\end{eqnarray*}
This system was obtained by considering two RLC sections
in Fig. \ref{fig:exam}(a) and the impedance matrix with
two ports at the end nodes.

In the SVD canonical form \rf{eqn:svd},
the rank was four and $I'_r = {\rm diag}(1, 1, -1, -1)$ was obtained.
The submatrix $A_{22}$ was given as follows:
\begin{eqnarray*}
\begin{aligned}
A_{22} &=
\left[
\begin{array}{rrr}
-2 & 1 & 0 \\
 1 & -1 & 0 \\
 0 & 0 & -1
\end{array}
\right].
\end{aligned}
\end{eqnarray*}
As this matrix is nonsingular, this system is index-1.
The SVD canonical form
was converted into state equation \rf{eqn:des11}
with the following coefficient matrices:
\begin{eqnarray*}
\begin{aligned}
\hat{A}_0 &=
\left[
\begin{array}{rrrr}
  0 & 0 & 1 & -1 \\
  0 & 0 & 0 & 1 \\
 -1 & 0 & -2 & 0 \\
  1 & -1 & 0 & -1
\end{array}
\right], \quad
\hat{B}_0 = -\hat{C}_0^T, \\
\hat{C}_0 &=
\left[
\begin{array}{rrrr}
  0 & 0 & -1 & 0 \\
  0 & 0 & -2 & 0
\end{array}
\right], \quad
\hat{D}_0 =
\left[
\begin{array}{rr}
  1 & 1 \\
  1 & 2
\end{array}
\right],
\end{aligned}
\end{eqnarray*}
where $I'_r \hat{A}_0$ is a symmetric matrix.
The AREs \rf{eqn:ric1} and \rf{eqn:ric2}
 were solved using the MATLAB {\tt care}() function, and
the following solutions were obtained:
\begin{eqnarray*}
\begin{aligned}
X &=
\left[
\begin{array}{rrrr}
  0.3439 & 0.1466 & -0.1298 & -0.1383 \\
  0.1466 & 0.2945 & 0.1298 & 0.0084 \\
 -0.1298 & 0.1298 & 0.4904 & 0.0804 \\
 -0.1383 & 0.0084 & 0.0804 & 0.1499
\end{array}
\right], \\
Y &=
\left[
\begin{array}{rrrr}
  0.3439 & 0.1466 & 0.1298 & 0.1383 \\
  0.1466 & 0.2945 & -0.1298 & -0.0084 \\
  0.1298 & -0.1298 & 0.4904 & 0.0804 \\
  0.1383 & -0.0084 & 0.0804 & 0.1499
\end{array}
\right].
\end{aligned}
\end{eqnarray*}
Here, $I'_r Y I'_r = X$; thus, Theorem \ref{th:are} holds.

After calculating the admittance and hybrid matrices,
we converted the admittance and hybrid matrices to the SVD canonical form.
Note that the rank was four for each case.
The submatrices $A_{22}$s for the admittance and hybrid matrices
were obtained respectively as:
\begin{eqnarray*}
\begin{aligned}
A_{22,y} &=
\left[
\begin{array}{rrrrr}
 -2 & 1 & 0 & -1 & 0 \\
  1 & -1 & 0 & 0 & 0 \\
  0 & 0 & -1 & 0 & 0 \\
 -1 & 0 & 0 & 0 & 0 \\
  0 & 0 & 0 & 0 & 0
\end{array}
\right], \quad
A_{22,h} =
\left[
\begin{array}{rrrr}
 -2 & 1 & 0 & 0 \\
  1 & -1 & 0 & 0 \\
  0 & 0 & -1 & 0 \\
  0 & 0 & 0 & 0
\end{array}
\right].
\end{aligned}
\end{eqnarray*}
Since these matrices are singular,
both systems are index-2; therefore, the index number depends on
the circuit structure and which parameter matrix is used.


\subsection{GARE and ARE Solutions for Index-2 Systems}

The second example was generated by considering two RLC sections
in Fig. \ref{fig:exam}(b) and the impedance matrix with
two ports at the end nodes.
As submatrix $A_{22}$ was singular,
the descriptor system \rf{eqn:des2} was obtained by
the spectral projectors \rf{eqn:proj}.
To define the GARE, $M_0$ of \rf{eqn:m0m1} was calculated as follows:
\begin{eqnarray*}
\begin{aligned}
M_0 &=
\left[
\begin{array}{rr}
 0 & 0 \\
 0 & 1
\end{array}
\right].
\end{aligned}
\end{eqnarray*}
As this is singular, matrix $D$ associated with
\rf{eqn:gric1} and \rf{eqn:gric2} was approximated by
$DD^T = \left( M_0 + M_0^T + \epsilon I_m \right)^{-1}$ with
$\epsilon = 1.0 \times 10^{-5}$.
Note that there is no MATLAB function to solve a GARE with singular matrix $E$;
thus, the solutions were computed using the Newton method \ct{benner3}.
From Theorem \ref{th:gare},
the GARE solution of \rf{eqn:gric1} is identical to that of \rf{eqn:gric2}.
The same solution $X=Y$ was obtained by the Newton method
as follows:
\begin{eqnarray}
X &=&
\left[
\begin{array}{rrrrr}
 1.0 \times 10^{-1} & 2.1 \times 10^{-4} & 1.7 \times 10^{-3} & 2.4 \times 10^{-20} & 2.1 \times 10^{-4} \\
 2.1 \times 10^{-4} & 4.0 \times 10^{-1} & -9.3 \times 10^{-2} & 6.8 \times 10^{-17} & 4.0 \times 10^{-1} \\
 1.7 \times 10^{-3} & -9.3 \times 10^{-2} & 2.4 \times 10^{-1} & -2.0 \times 10^{-17} & -9.3 \times 10^{-2} \\
 2.4 \times 10^{-20} & 6.8 \times 10^{-17} & -2.0 \times 10^{-17} & 1.0 \times 10^{-32} & 6.8 \times 10^{-17} \\
 2.1 \times 10^{-4} & 4.0 \times 10^{-1} & -9.3 \times 10^{-2} & 6.8 \times 10^{-17} & 4.0 \times 10^{-1}
\end{array}
\right] \nonumber \\
&& \hspace*{-5mm}
+ j
\left[
\begin{array}{rrrrr}
 9.5 \times 10^{-17} & 5.6 \times 10^{-17} & 9.7 \times 10^{-17} & 9.2 \times 10^{-33} & 1.1 \times 10^{-16} \\
-1.1 \times 10^{-16} & -1.1 \times 10^{-13} & 4.2 \times 10^{-14} & 4.6 \times 10^{-29} & 2.7 \times 10^{-13} \\
-1.4 \times 10^{-17} & 9.5 \times 10^{-15} & -2.4 \times 10^{-14} & 3.5 \times 10^{-29} & 2.6 \times 10^{-13} \\
-2.1 \times 10^{-24} & 2.4 \times 10^{-21} & 9.5 \times 10^{-22} & 2.7 \times 10^{-37} & 2.4 \times 10^{-21} \\
-1.9 \times 10^{-16} & -1.0 \times 10^{-13} & 4.2 \times 10^{-14} & 4.6 \times 10^{-29} & 2.7 \times 10^{-13}
\end{array}
\right], \nonumber \\
\label{eqn:ex2_g}
\end{eqnarray}
where the Lyapunov equation at each Newton step
was solved by the ADI method \ct{stykel4}, and the eigenvalues obtained
by full decomposition were used as the shift parameters of the ADI method.
Here, the imaginary part of the solution can be ignored;
thus, the matrix is considered symmetric.
The solutions for the admittance and hybrid matrices were also calculated,
and the solutions were obtained as $X=Y$.

A realization for an index-2 system is obtained in Section 4.2,
where the above example with the impedance matrix is used.
In this case, the transfer matrix is expressed as \rf{eqn:elm},
and matrix $I'_r A_s$ is obtained as follows:
\begin{eqnarray*}
I'_r A_s &=&
\left[
\begin{array}{rrrr}
-1 & 0 & -1 & 0 \\
 0 & 0 & 1 & 0 \\
 1 & -1 & -1 & 0 \\
 0 & 0 & 0 & 0
\end{array}
\right].
\end{eqnarray*}
As the rank of the matrix is three, no value of the
transfer function at $s=0$ exists, which contradicts the original system.
Fortunately, as $\hat{C}_2 = \hat{B}_2^T = 0$, the second term
of the last line of \rf{eqn:elm} can be ignored; therefore,
the transfer function does have a value at $s=0$.

After obtaining the AREs, we solved them using the MATLAB
{\tt care}() function.
The solutions of \rf{eqn:ric1} and \rf{eqn:ric2} were
obtained respectively as follows:
\begin{eqnarray}
\begin{aligned}
X &=
\left[
\begin{array}{rrr}
  9.96 \times 10^{-1} & 2.12 \times 10^{-4} & -1.69 \times 10^{-3} \\
  2.12 \times 10^{-4} & 3.97 \times 10^{-1} & 9.33 \times 10^{-2} \\
 -1.69 \times 10^{-3} & 9.33 \times 10^{-2} & 2.39 \times 10^{-1}
\end{array}
\right], \\
Y &=
\left[
\begin{array}{rrr}
 9.96 \times 10^{-1} & 2.12 \times 10^{-4} & 1.69 \times 10^{-3} \\
 2.12 \times 10^{-4} & 3.97 \times 10^{-1} & -9.33 \times 10^{-2} \\
 1.69 \times 10^{-3} & -9.33 \times 10^{-2} & 2.39 \times 10^{-1}
\end{array}
\right].
\end{aligned}
\label{eqn:ex2_a}
\end{eqnarray}
Here, $Y = I'_r X I'_r$ (precisely $Y = I''_r X I''_r$) is confirmed and
Theorem \ref{th:are} holds.
The eigenvalues of the real part of \rf{eqn:ex2_g}
are
$9.96 \times 10^{-1}$, $2.09 \times 10^{-1}$, $8.23 \times 10^{-1}$,
$-3.60 \times 10^{-13}$, and $-1.41 \times 10^{-33}$,
and those of \rf{eqn:ex2_a} are $1.96 \times 10^{-1}$,
$4.40 \times 10^{-1}$, and $9.96 \times 10^{-1}$.
The fourth and fifth largest eigenvalues of \rf{eqn:ex2_g}
are considered to be zero; thus,
nonnegative solutions were obtained.


\subsection{Frequency Response Errors and Residual Errors of ARE and GARE}

\begin{figure}[tb]
\begin{center}
\subfigure[]{\includegraphics[width=0.43\columnwidth]{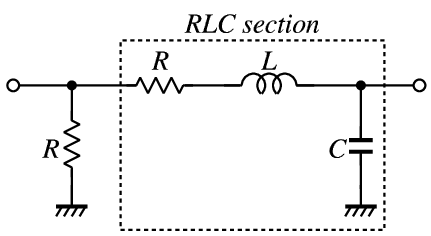}}
\subfigure[]{\includegraphics[width=0.43\columnwidth]{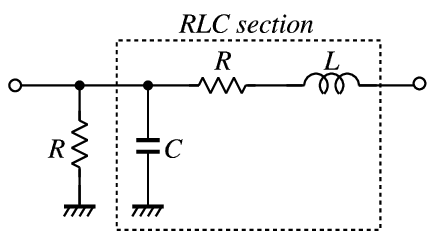}}
\end{center}
\caption{Fundamental circuit structures}
\label{fig:exam}
\end{figure}

Two examples were analyzed to evaluate the proposed method.
The fundamental circuit structures are shown in Figs. \ref{fig:exam}(a)
and \ref{fig:exam}(b), where $R = 1[\Omega]$, $L = 1[\mbox{nH}]$,
and $L = 1[\mbox{nF}]$.
Note that 100 RLC sections were considered in the numerical examples.
The first example was obtained by the voltage-current relationship of
the leftmost two nodes of the circuit with 100 RLC sections in Fig. \ref{fig:exam}(a);
thus, $m = 2$ in \rf{eqn:des1}. Here, this system has index-1.
The second example was obtained by the relationship of the left end of the
first RLC section in Fig. \ref{fig:exam}(b) and the right end of the second
RLC section; thus, $m = 2$. Here, this system has index-2.
We refer to these as index-1 and index-2 examples, respectively.

\begin{figure}[tb]
\begin{center}
\subfigure[]{\includegraphics[width=0.43\columnwidth]{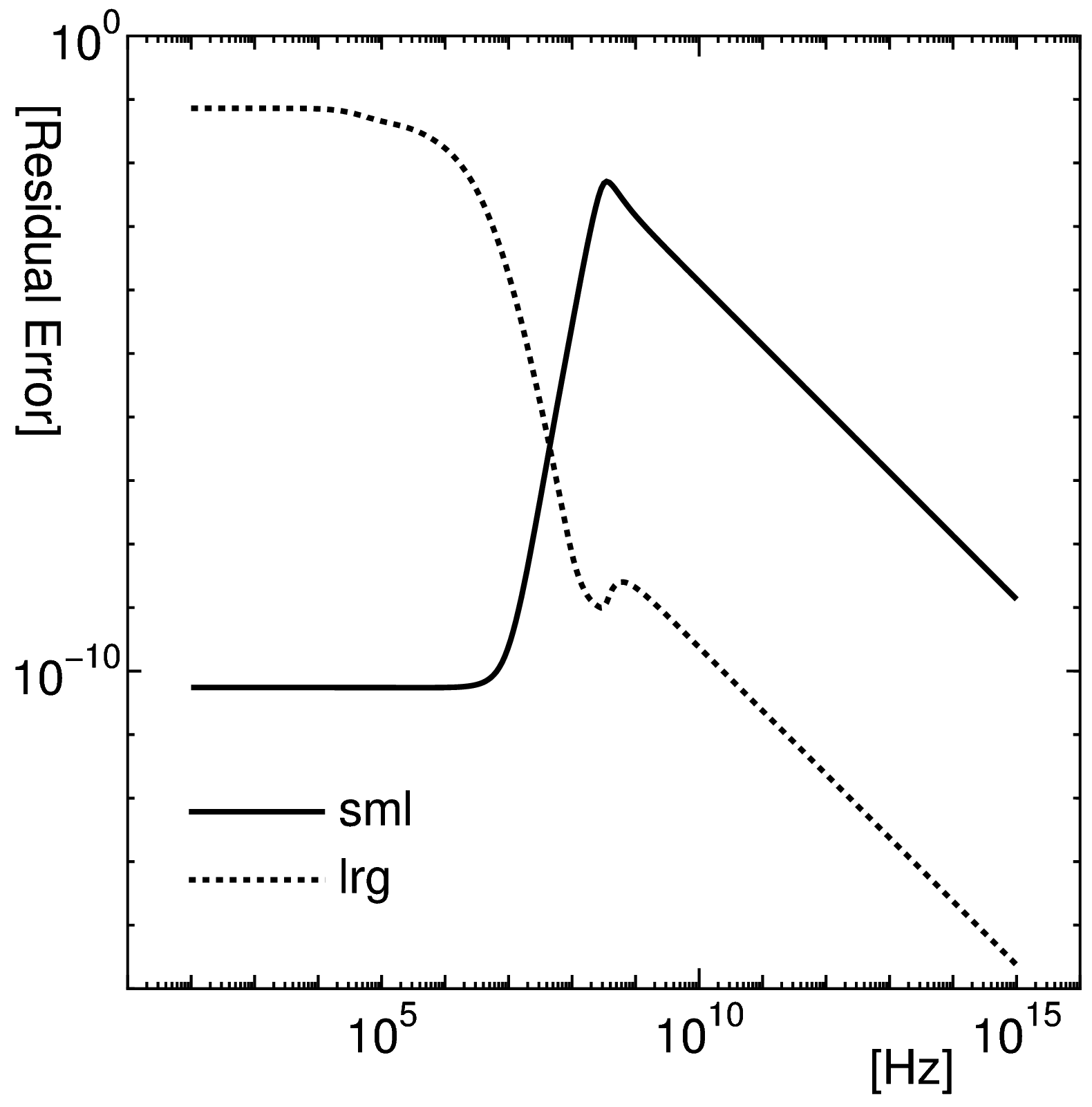}}
\subfigure[]{\includegraphics[width=0.43\columnwidth]{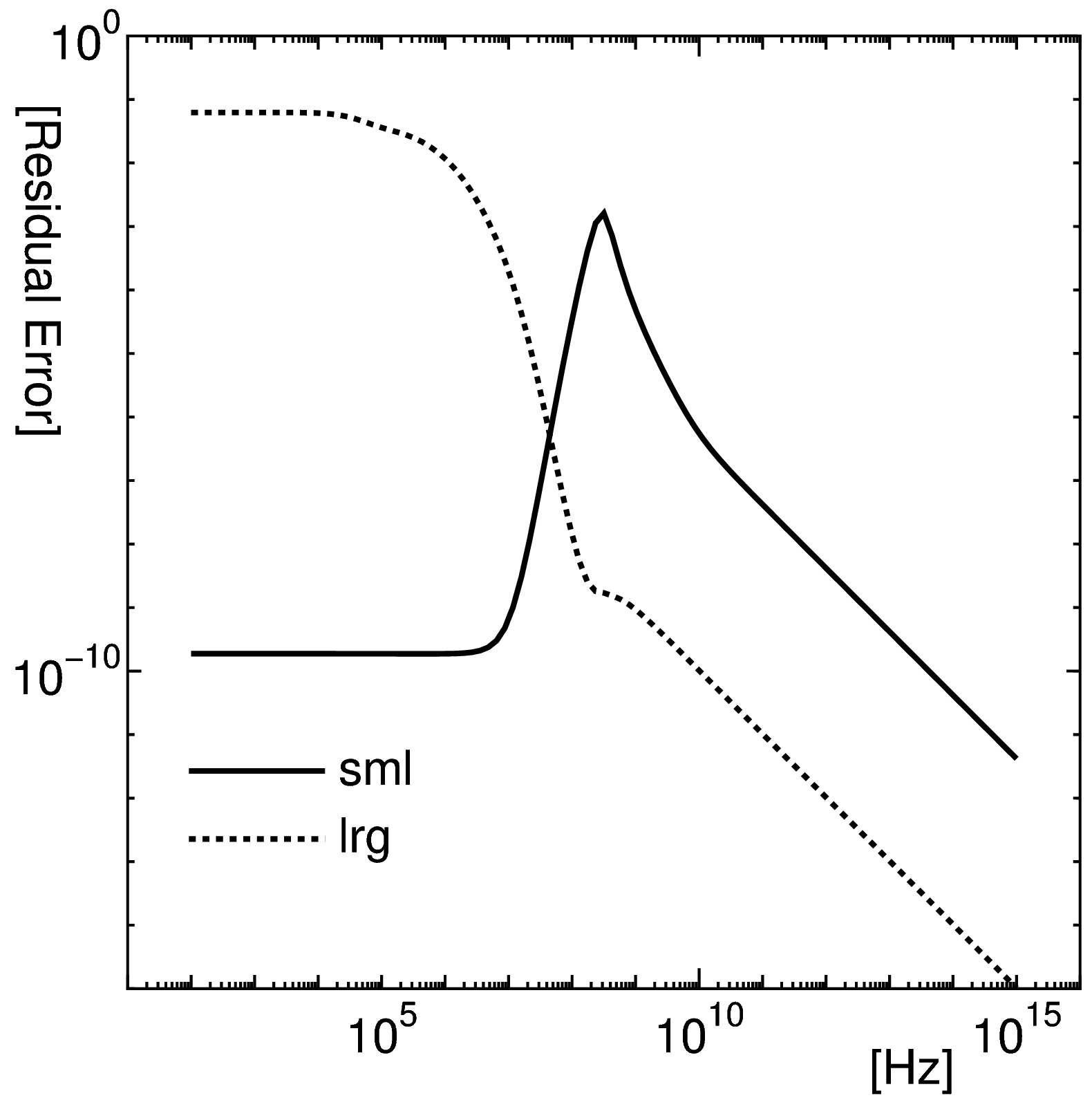}}
\end{center}
\caption{Relative frequency response errors of
reduced-order models for index-1 example obtained by (a){\bf SPRBT-1} and
(b){\bf SPRBT-2}}
\label{fig:error1}
\end{figure}

\begin{figure}[tb]
\begin{center}
\subfigure[]{\includegraphics[width=0.43\columnwidth]{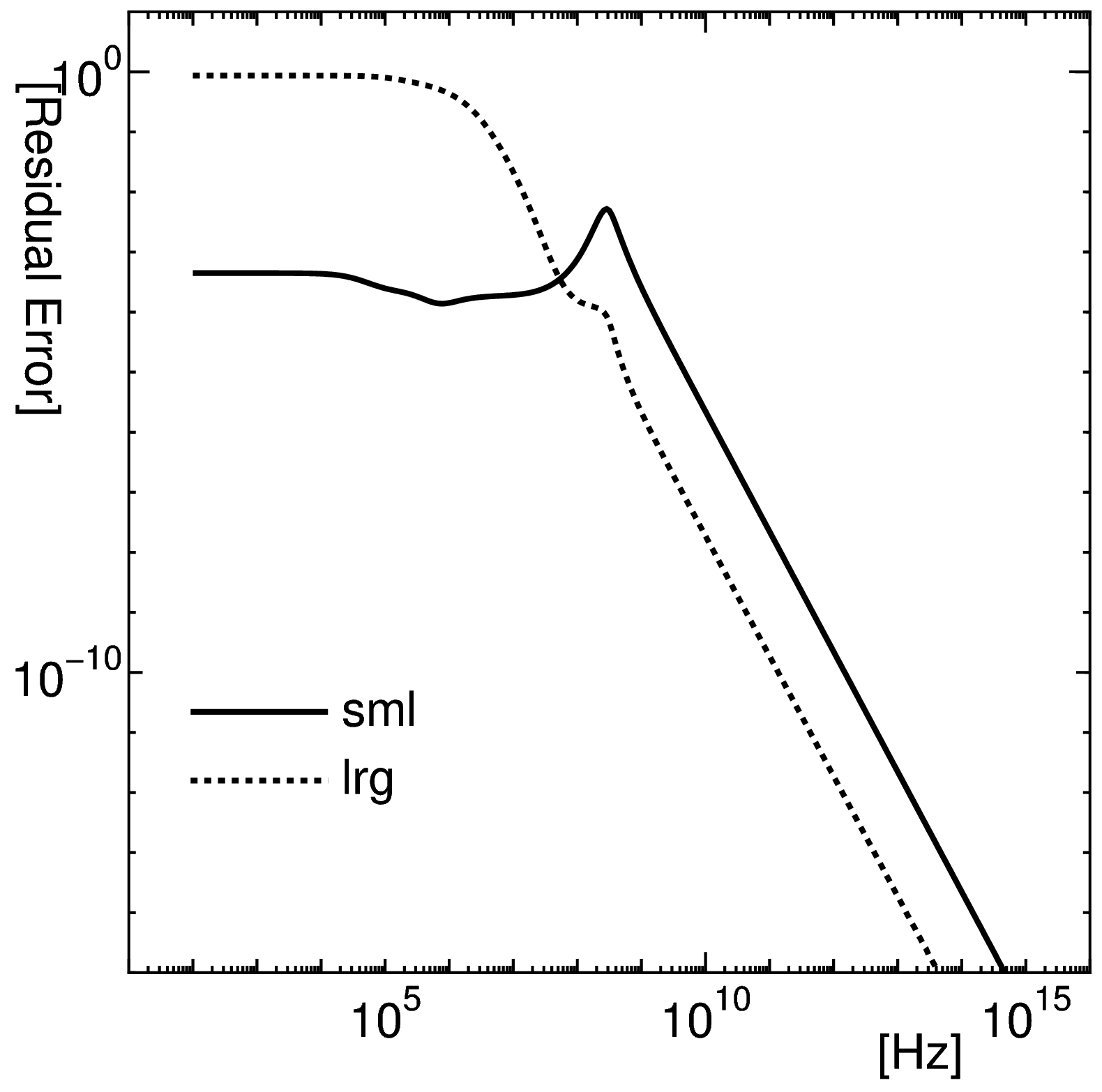}}
\subfigure[]{\includegraphics[width=0.43\columnwidth]{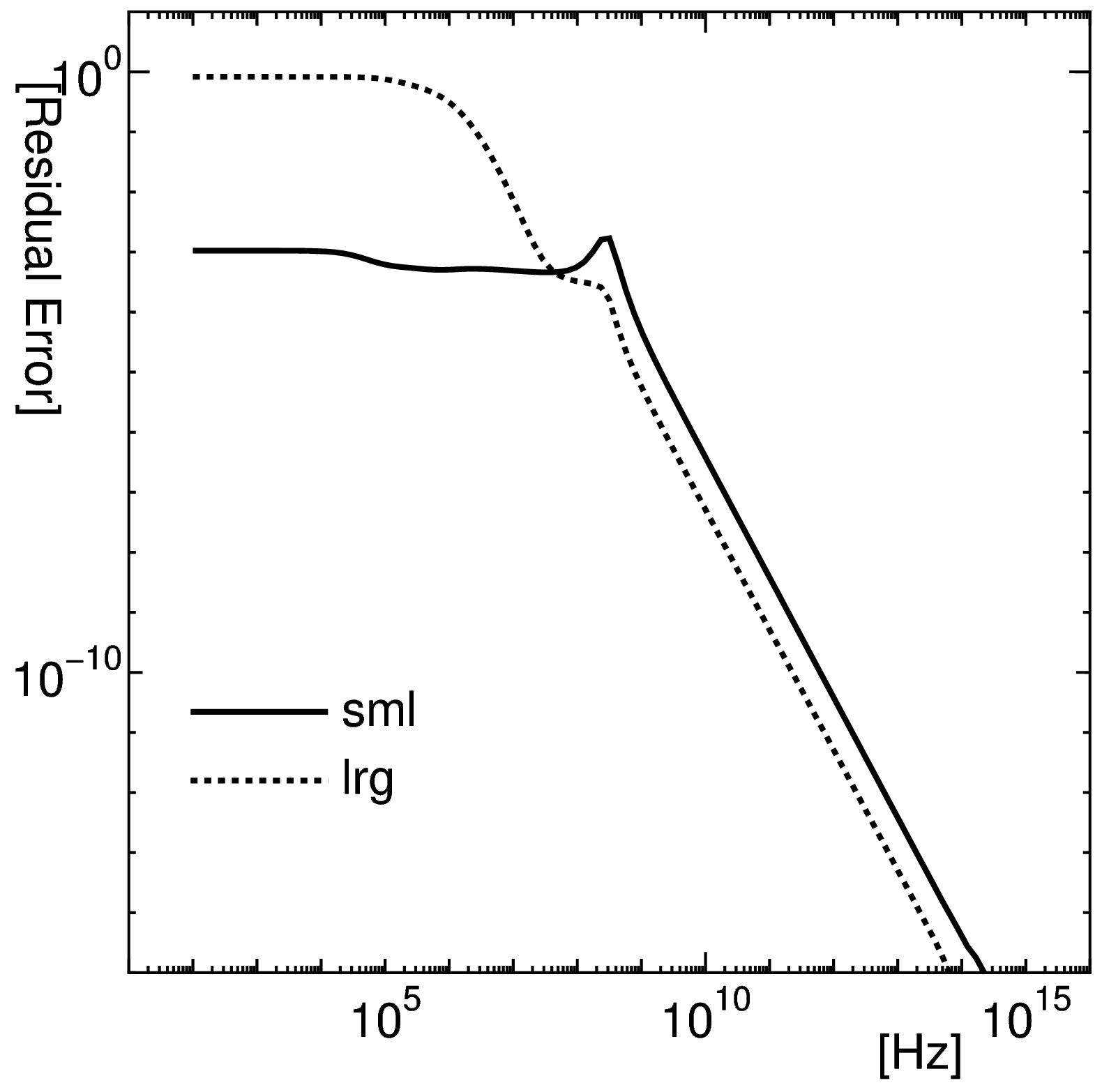}}
\end{center}
\caption{Relative frequency response errors of
reduced-order models for index-2 example obtained by (a){\bf SPRBT-1} and
(b){\bf SPRBT-2}}
\label{fig:error2}
\end{figure}

We calculated $30$ eigenvalues of
the inverse Hamiltonian matrix ${\cal H}^{-1}$ for {\bf SPRBT-1}
and ${\cal H}^{-1} {\cal E}$ for {\bf SPRBT-2}
to obtain the shifts of RADI.
Note that where reciprocal values were used, we use the label "sml."
For comparison, $30$ eigenvalues of the Hamiltonian matrix ${\cal H}$
were also calculated using the Arnoldi method for {\bf SPRBT-1},
in which the obtained values
correspond to the large eigenvalues of the Hamiltonian matrix ${\cal H}$.
For {\bf SPRBT-2}, the shifts were calculated as follows \ct{tanji}.
Here, the generalized eigenvalue problem is expressed as
${\cal E} x = (1/\lambda) {\cal H} x = t {\cal H} x$.
Assuming an expansion point at $s = s_0$, the problem can be rewritten as
$({\cal E} - s_0 {\cal H}) x = (t-s_0) {\cal H} x$.
The matrix pencil $s {\cal E} - {\cal H}$ is regular;
thus, ${\cal E}-s_0 {\cal H}$ is nonsingular.
Therefore, the problem is described by
$({\cal E}-s_0 {\cal H})^{-1} {\cal H} x = 1/(t-s_0) x
= \xi x$.
Using the Arnoldi method, we obtain the approximate eigenvalues as
$\lambda = 1/(1/\xi + s_0)$.
Using the Krylov subspace method, large eigenvalues $\xi$ of matrix
$({\cal E}-s_0 {\cal H})^{-1} {\cal H}$
are obtained; thus, $\lambda$ also corresponds to a large generalized eigenvalue
of $({\cal E}, {\cal H})$.
As $t$ should be small, $s_0$ is assumed to be a small negative value.
Note that where these shifts were used, we use the label "lrg."

Figures \ref{fig:error1}(a) and \ref{fig:error1}(b) shows the frequency response
errors obtained by {\bf SPRBT-1} and {\bf SPRBT-2}, respectively, for the
index-1 example, in which 30 RADI steps with 15 shifts were applied and
$15$ order models were generated.
The responses of sml are satisfactory for both {\bf SPRBT-1} and {\bf SPRBT-2},
which implies that the suitable shifts were obtained.
The responses of lrg with $s_0 = -10^{-5}$ are inaccurate at low frequencies;
however,
they are accurate at high frequencies, which implies that
large eigenvalues contribute to model accuracy at high frequencies.
Figures \ref{fig:error2}(a) and \ref{fig:error2}(b) show the frequency response
errors obtained by {\bf SPRBT-1} and {\bf SPRBT-2}, respectively, for the
index-2 example, in which 30 RADI steps with 15 shifts were applied and
$15$ order models were generated.
Here, the responses of sml are accurate at low frequencies, which implies that
shift selection with Algorithm 4 is also suitable for index-2
descriptor systems.


\begin{table}[bt]
\caption{Computational time and memory usage of {\bf SPRBT-1} with (a) RARI
and (b) QADI for \underline{ index-1 example}}
{\footnotesize
\begin{center}
\subtable[]{
\begin{tabular}{|r|r|r|r|r|r|} \hline
  size & ADI [s] & Else [s] & Total [s] & memory [MB] \\ \hline
  $3,001$ & $0.11$ & $0.68$ & $0.79$ & $1.50 \times 10^0$ \\ \hline
  $300,001$ & $12.25$ & $36.44$ & $48.69$ & $1.50 \times 10^2$ \\ \hline
$30,000,001$ & $1,289.85$ & $3,657.96$ & $4,947.81$ & $1.50 \times 10^4$ \\ \hline
\end{tabular}}
\subtable[]{
\begin{tabular}{|r|r|r|r|r|r|} \hline
  size & ADI [s] & Else [s] & Total [s] & memory [MB] \\ \hline
  $3,001$ & $0.51$ & $0.57$ & $1.08$ & $1.50 \times 10^0$ \\ \hline
  $300,001$ & $48.41$ & $36.54$ & $84.95$ & $1.50 \times 10^2$ \\ \hline
$30,000,001$ & $5,399.45$ & $3,629.05$ & $9,028.50$ & $1.50 \times 10^4$ \\ \hline
\end{tabular}}
\end{center}}
\label{tab:id1_are}
\end{table}

\begin{table}[bt]
\caption{Computational time and memory usage of {\bf SPRBT-2} with (a) RARI
and (b) QADI for \underline{ index-1 example}}
{\footnotesize
\begin{center}
\subtable[]{
\begin{tabular}{|r|r|r|r|r|r|} \hline
  size & ADI [s] & Else [s] & Total [s] & memory [MB] \\ \hline
  $3,001$ & $0.11$ & $0.62$ & $0.73$ & $1.60 \times 10^0$ \\ \hline
  $300,001$ & $11.66$ & $39.08$ & $50.74$ & $1.60 \times 10^2$ \\ \hline
$30,000,001$ & $1,202.00$ & $3,721.13$ & $4,923.13$ & $1.59 \times 10^4$ \\ \hline
\end{tabular}}
\subtable[]{
\begin{tabular}{|r|r|r|r|r|r|} \hline
  size & ADI [s] & Else [s] & Total [s] & memory [MB] \\ \hline
  $3,001$ & $0.66$ & $0.62$ & $1.28$ & $1.60 \times 10^0$ \\ \hline
  $300,001$ & $71.81$ & $40.50$ & $112.31$ & $1.60 \times 10^2$ \\ \hline
$30,000,001$ & $7,838.62$ & $3,889.46$ & $11,728.07$ & $1.59 \times 10^4$ \\ \hline
\end{tabular}}
\end{center}}
\label{tab:id1_gare}
\end{table}


\begin{table}[bt]
\caption{Computational time and memory usage of {\bf SPRBT-1} with (a) RARI
and (b) QADI for \underline{ index-2 example}}
{\footnotesize
\begin{center}
\subtable[]{
\begin{tabular}{|r|r|r|r|r|r|} \hline
  size & ADI [s] & Else [s] & Total [s] & memory [MB] \\ \hline
  $3,001$ & $0.12$ & $0.46$ & $0.57$ & $2.74 \times 10^0$ \\ \hline
  $300,001$ & $14.09$ & $34.73$ & $48.82$ & $2.74 \times 10^2$ \\ \hline
$30,000,001$ & $2,326.66$ & $5,134.52$ & $7,461.18$ & $2.74 \times 10^4$ \\ \hline
\end{tabular}}
\subtable[]{
\begin{tabular}{|r|r|r|r|r|r|} \hline
  size & ADI [s] & Else [s] & Total [s] & memory [MB] \\ \hline
  $3,001$ & $0.63$ & $0.51$ & $1.14$ & $2.70 \times 10^0$ \\ \hline
  $300,001$ & $59.13$ & $34.72$ & $93.85$ & $2.74 \times 10^2$ \\ \hline
$30,000,001$ & $18,073.48$ & $4,984.57$ & $23,058.05$ & $2.74 \times 10^4$ \\ \hline
\end{tabular}}
\end{center}}
\label{tab:id2_are}
\end{table}

\begin{table}[bt]
\caption{Computational time and memory usage of {\bf SPRBT-2} with (a) RARI
and (b) QADI for \underline{ index-2 example}}
{\footnotesize
\begin{center}
\subtable[]{
\begin{tabular}{|r|r|r|r|r|r|} \hline
  size & ADI [s] & Else [s] & Total [s] & memory [MB] \\ \hline
  $3,001$ & $0.16$ & $0.79$ & $0.95$ & $2.55 \times 10^0$ \\ \hline
  $300,001$ & $18.39$ & $62.71$ & $81.11$ & $2.54 \times 10^2$ \\ \hline
$30,000,001$ & $3,005.42$ & $7,608.21$ & $10,613.63$ & $2.54 \times 10^4$ \\ \hline
\end{tabular}}
\subtable[]{
\begin{tabular}{|r|r|r|r|r|r|} \hline
  size & ADI [s] & Else [s] & Total [s] & memory [MB] \\ \hline
  $3,001$ & $0.73$ & $0.83$ & $1.56$ & $2.55 \times 10^0$ \\ \hline
  $300,001$ & $76.73$ & $60.06$ & $136.79$ & $2.54 \times 10^2$ \\ \hline
$30,000,001$ & $19,702.38$ & $7,3497.40$& $27,051.78$ & $2.54 \times 10^4$ \\ \hline
\end{tabular}}
\end{center}}
\label{tab:id2_gare}
\end{table}


\subsection{Computational Time and Memory Usage}

For the index-1 and index-2 examples, we measured the
computational time and memory usage when applying {\bf RPRBT-1}
and {\bf RPRBT-2}.
The simulations were performed on a computer with a 3.7-GHz
Intel Xeon E5-1620 CPU and 32 GB of memory.
In Tables \ref{tab:id1_are}-\ref{tab:id2_gare},
"size" is the order of the descriptor system \rf{eqn:des1},
"ADI" is the computational time of RADI or QADI,
"Else" is the time except for RADI or QADI,
"Total" is the total calculation time,
and "Mem" is memory usage.

{\bf RPRBT-1} and {\bf RPRBT-2} were applied to the index-1 example,
in which 20 RADI steps with 10 shifts were applied and 15 order models were
generated.
Tables \ref{tab:id1_are} and \ref{tab:id1_gare} show the computational time
and memory usage for {\bf RPRBT-1} and {\bf RPRBT-2}, respectively,
in which RADI is compared to QADI.
As shown in Table \ref{tab:id1_are},
RADI is 4.2 times faster than QADI, and the total time using RADI
is 1.8 times less than when using QADI.
As shown in Table \ref{tab:id1_gare},
RADI is 6.52 times faster than QADI, and the total time using RADI
is 2.4 times less than when using QADI.

{\bf RPRBT-1} and {\bf RPRBT-2} were also applied to the index-2 example,
in which 20 RADI steps with 10 shifts were applied and 15 order models were
generated.
Tables \ref{tab:id2_are} and \ref{tab:id2_gare} show the computational time
and memory usage for {\bf RPRBT-1} and {\bf RPRBT-2}, respectively.
As shown in Table \ref{tab:id2_are},
RADI is 7.7 times faster than QADI, and the total time using RADI
is 3.1 times less than when using QADI.
As shown in Table \ref{tab:id2_gare},
RADI is 6.5 times faster than QADI, and the total time using RADI
is 2.6 times less than when using QADI.

For the index-1 example, {\bf RPRBT-1} is nearly identical to {\bf RPRBT-2}
relative to efficiency and memory usage when RADI was used.
On the other hand, for the index-2 example,
{\bf RPRBT-1} is 1.4 times faster than {\bf RPRBT-2} when RADI was used, and
the memory usage of {\bf RPRBT-1} is compatible with that of {\bf RPRBT-2}.

Consequently, solving AREs or GAREs is no longer dominant
in the computational cost of PRBT.
Moreover, the efficiency and memory usage of {\bf RPRBT-1} are nearly compatible
to {\bf RPRBT-2} when RADI was used.

%
%

\section{Conclusions}

In this paper, we have presented reciprocal and PRBTs
for index-1 and index-2 descriptor systems,
for which two approaches based on AREs or GAREs have been proposed.
Furthermore, RADI was introduced to solve AREs and GAREs.
We have demonstrated that solving the Riccati equation is
no longer dominant in the computational cost of
PRBT.
By comparing the two approaches, we have also shown that both methods
are compatible relative to computational time and memory usage.
In addition, some properties of ARE- and
GARE-associated descriptor systems for passive electrical circuits
have been provided.

%
%

%
%

\section*{Appendix A}

Here, we prove Theorem \ref{th:are}.
From the symmetry of $A_0$, $A_{11} = A_{11}^T$,
$A_{22} = A_{22}^T$, and $A_{21} = A_{12}^T$,
we can obtain $\hat{A}_0 = I'_r A_s$, where $A_s$ is a symmetric matrix.
From \rf{eqn:sec2}, we can obtain $\hat{B}_0 = I'_r \hat{C}_0^T$
for an impedance matrix and
$\hat{B}_0 = -I'_r \hat{C}_0^T$ for an admittance matrix.
Then, the AREs \rf{eqn:ric1} and \rf{eqn:ric2} are
rewritten respectively as follows:
\begin{eqnarray}
&& \hspace*{-5mm} \left( I'_r A_s \mp I'_r C_0^T D D^T C_0 \right)^T X
+ X \left( I'_r A_s \mp I'_r C_0^T D D^T C_0 \right) \nonumber \\
&& + X I'_r C_0^T D D^T C_0 I'_r X + C_0^T D D^T C_0 = 0,
\label{eqn:ric1_2} \\
&& \hspace*{-5mm} \left( I'_r A_s \mp I'_r C_0^T D D^T C_0 \right)^T
I'_r Y I'_r
+ I'_r Y I'_r \left( I'_r A_s \mp I'_r C_0^T D D^T C_0 \right) \nonumber \\
&& + I'_r Y I'_r I'_r C_0^T D D^T C_0 I'_r I'_r Y I'_r + C_0^T D D^T C_0 = 0,
\label{eqn:ric2_2}
\end{eqnarray}
where the minus and plus signs correspond to
the impedance and admittance matrices, respectively.
From \rf{eqn:ric1_2} and \rf{eqn:ric2_2}, $X = I'_r Y I'_r$, which indicates
that $Y = I'_r X I'_r$.

For a hybrid matrix,
we have $\hat{B}_0 = I'_r \hat{C}_0^T I^\circ_m$ from
\rf{eqn:gh} and \rf{eqn:sec2}.
As $\hat{D}_0 = - I^\circ_m \hat{B}_2^T A_{22}^{-1} \hat{B}_2$,
$D D^T = ( \hat{D}_0 + \hat{D}_0^T )^{-1}$ becomes a block
diagonal matrix; thus, the relationship $I^\circ_m D D^T I^\circ_m = D D^T$
is obtained.
Using this relationship, the dual AREs are expressed as follows:
\begin{eqnarray}
&& \hspace*{-5mm} \left( I'_r A_s - I'_r C_0^T I^\circ_m D D^T C_0 \right)^T X
+ X \left( I'_r A_s - I'_r C_0^T I^\circ_m D D^T C_0 \right) \nonumber \\
&& + X I'_r C_0^T D D^T C_0 I'_r X + C_0^T D D^T C_0 = 0,
\label{eqn:ric1_3} \\
&& \hspace*{-5mm} \left( I'_r A_s - I'_r C_0^T I^\circ_m D D^T C_0 \right)^T I'_r Y I'_r
+ I'_r Y I'_r \left( I'_r A_s - I'_r C_0^T I^\circ_m D D^T C_0 \right) \nonumber \\
&& + I'_r Y I'_r I'_r C_0^T D D^T C_0 I'_r I'_r X I'_r + C_0^T D D^T C_0 = 0.
\label{eqn:ric2_3}
\end{eqnarray}
From \rf{eqn:ric1_3} and \rf{eqn:ric2_3}, we obtain $X = I'_r Y I'_r$,
which indicates that $Y = I'_r X I'_r$.


\section*{Appendix B}

Here, we prove the equations in \rf{eqn:m0m1}.
Using \rf{eqn:weire}, we obtain the following relation:
\begin{eqnarray*}
&& \hspace*{-10mm}
- \bar{C}_0 \left( I_n - P_r \right) \bar{A}_0^{-1}
\left( I_n - P_l \right) \bar{B}_0 \\
&=& - \bar{C}_0 T_r^{-1}
\left[
\begin{array}{cc}
0 & 0 \\
0 & I_{n-r}
\end{array}
\right]
\left[
\begin{array}{cc}
J^{-1} & 0 \\
0 & I_{n-r}
\end{array}
\right]
\left[
\begin{array}{cc}
0 & 0 \\
0 & I_{n-r}
\end{array}
\right]
T_l^{-1} \bar{B}_0 \\
&=& - C_{\infty} B_{\infty}.
\end{eqnarray*}
Then, $M_0$ is obtained by adding $\bar{D}_0$ to it.
Similarly, $M_1$ of \rf{eqn:m0m1} is given as follows:
\begin{eqnarray*}
&& \hspace*{-10mm}
- \bar{C}_0 \bar{A}_0^{-1} \left( I_n - P_l \right) \bar{E}_0
\left( I_n - P_r \right) \bar{A}_0^{-1} \bar{B}_0 \\
&=& - \bar{C}_0 T_r^{-1}
\left[
\begin{array}{cc}
J^{-1} & 0 \\
0 & I_{n-r}
\end{array}
\right]
\left[
\begin{array}{cc}
0 & 0 \\
0 & I_{n-r}
\end{array}
\right] \\
&& \times
\left[
\begin{array}{cc}
I_r & 0 \\
0 & N
\end{array}
\right]
\left[
\begin{array}{cc}
0 & 0 \\
0 & I_{n-r}
\end{array}
\right]
\left[
\begin{array}{cc}
J^{-1} & 0 \\
0 & I_{n-r}
\end{array}
\right]
T_l^{-1} \bar{B}_0 \\
&=& - C_{\infty} N B_{\infty}.
\end{eqnarray*}
%


\section*{Appendix C}
Here, we prove Theorem \ref{th:gare} with the dual projected generalized
Lur'e equations rather than dual GAREs \rf{eqn:gric1}
and \rf{eqn:gric2}.
The dual GAREs \rf{eqn:gric1} and \rf{eqn:gric2} are respectively
equivalent to the projected generalized Lur'e equations:
\begin{eqnarray}
\left\{
\begin{aligned}
\bar{A}_0 ^T X \bar{E}_0 + \bar{E}_0 ^T X \bar{A}_0 &= -K_o^T K_o, \quad X = P_l^T X P_l \\
\bar{E}_0 ^T X \bar{B}_0 - P_r^T \bar{C}_0^T &= -K_o^T J_o, \quad M_o + M_o^T = J_o^T J_o
\end{aligned}
\right.
\label{eqn:lur1} \\
\left\{
\begin{aligned}
\bar{A}_0 Y \bar{E}_0^T + \bar{E}_0 Y \bar{A}_0^T &= -K_c K_c^T, \quad Y = P_r Y P_r^T \\
\bar{E}_0 Y \bar{C}_0^T - P_l \bar{B}_0 &= -K_c J_c^T, \quad M_o + M_o^T = J_c J_c^T
\end{aligned}
\right.
\label{eqn:lur2}
\end{eqnarray}

For an impedance matrix, $\bar{E}_0=\bar{E}_0^T$, $\bar{A}_0=\bar{A}_0^T$,
and $\bar{B}_0=\bar{C}_0^T$ hold.
Using the relationship $P_l = P_r^T$ and rewriting $K_c$ and $J_c$ with
$K_o^T$ and $J_o^T$ in \rf{eqn:lur2}, respectively, we confirm that $X = Y$.
Similar to the impedance matrix case, $X = Y$ for an admittance matrix.

For a hybrid matrix, $\bar{E}_0=\bar{E}_0^T$, $\bar{A}_0=\bar{A}_0^T$,
and $\bar{B}_0=\bar{C}_0^T I^\circ_m$
hold. Using the relationship $P_l = P_r^T$,
we can write the third equation of \rf{eqn:lur1} as
$\bar{E}_0 X \bar{C}_0^T I^\circ_m - P_l \bar{B}_0 I^\circ_m = -K^T_o J_o$.
Here, $M_0$ is block skew symmetric; thus, $J_o^T J_o = I^\circ_m J_o^T J_o I^\circ_m$.
This implies that $J_o = J_o I^\circ_m$.
Therefore, the dual Lur'e equations are identical
and $X = Y$ for the hybrid matrix.

\end{document}